\documentclass[11pt,reqno]{amsart}

\usepackage{amsmath, amsfonts, amssymb, amscd, enumerate, young, youngtab, empheq, stmaryrd, slashed,
graphicx, mathrsfs, graphpap, curves, color, mathabx}
\theoremstyle{plain}
\newtheorem{theo}{Theorem}[section]

\theoremstyle{definition}

\newtheorem{definition}[theo]{Definition}

\theoremstyle{plain}
\newtheorem{lemma}[theo]{Lemma}
\newtheorem{theorem}[theo]{Theorem}
\newtheorem{corollary}[theo]{Corollary}
\newtheorem{proposition}[theo]{Proposition}

\theoremstyle{definition}

\newtheorem{remark}[theo]{Remark}

\newcommand{\rank}{\operatorname{rank}}

\newcommand{\beq}{\begin{equation}}
\newcommand{\eeq}{\end{equation}}
\renewcommand{\a}{\alpha}

\renewcommand{\o}{\omega}

\newcommand{\tE}{\textbf{E}}
\newcommand{\tH}{\textbf{H}}

\renewcommand{\O}{\Omega}

\newcommand{\bC}{\mathbb{C}}

\newcommand{\bR}{\mathbb{R}}

\newcommand{\bH}{\mathbb{H}}

\newcommand{\tD}{\textbf{D}}

\newcommand{\tgh}{\wt{\mathfrak{h}}}
\newcommand{\tgk}{\wt{\mathfrak{k}}}

\newcommand{\gc}{\mathfrak{c}}

\renewcommand{\gg}{\mathfrak{g}}
\newcommand{\gh}{\mathfrak{h}}
\newcommand{\gk}{\mathfrak{k}}

\newcommand{\gm}{\mathfrak{m}}

\newcommand{\gs}{\mathfrak{s}}

\newcommand{\so}{\mathfrak{so}}

\newcommand{\ggl}{\mathfrak{gl}}
\declareslashed{}{/}{0}{0}{\partial}

\makeatletter
\def\widebreve#1{\mathop{\vbox{\m@th\ialign{##\crcr\noalign{\kern\p@}%
  \brevefill\crcr\noalign{\kern0.1\p@\nointerlineskip}%
  $\hfil\displaystyle{#1}\hfil$\crcr}}}\limits}
\def\brevefill{$\m@th \setbox\z@\hbox{}%
 \hfill\scalebox{0.6}{\rotatebox[origin=c]{90}{(}} \kern0pt $}
\makeatletter

\newcommand\GL{\mathrm{GL}}
\newcommand\SL{\mathrm{SL}}
\newcommand\SO{\mathrm{SO}}
\newcommand\CO{\mathrm{CO}}
\newcommand\OO{\mathrm{O}}

\newcommand\Sp{\mathrm{Sp}}
\renewcommand\sp{\mathfrak{sp}}
\renewcommand\sl{\mathfrak{sl}}

\newcommand{\omc}{\varpi}
\newcommand{\cA}{\mathcal{A}}

\newcommand{\cC}{\mathcal{C}}
\newcommand{\cD}{\mathcal{D}}

\newcommand{\cR}{\mathcal{R}}

\newcommand{\cT}{\mathcal{T}}
\newcommand{\cU}{\mathcal{U}}

\renewcommand{\square}{\kern1pt\vbox
{\hrule height 0.6pt\hbox{\vrule width 0.6pt\hskip 3pt
\vbox{\vskip 6pt}\hskip 3pt\vrule width 0.6pt}\hrule height0.6pt}\kern1pt}

\DeclareMathOperator\End{End\;}

\DeclareMathOperator\Ad{Ad}
\DeclareMathOperator\ad{ad}

\newcommand{\wt}{\widetilde}
\newcommand{\wh}{\widehat}

\newcommand{\be}{\begin{equation}}
\newcommand{\ee}{\end{equation}}

\def\<#1,#2>{\langle\,#1,\,#2\,\rangle}
\newcommand{\arr}{\begin{array}{rlll}}
\newcommand{\ea}{\end{array}}
\newcommand{\bea}{\begin{eqnarray}}
\newcommand{\eea}{\end{eqnarray}}
\newcommand{\bean}{\begin{eqnarray*}}
\newcommand{\eean}{\end{eqnarray*}}

\def\sideremark#1{\ifvmode\leavevmode\fi\vadjust{
\vbox to0pt{\hbox to 0pt{\hskip\hsize\hskip1em
\vbox{\hsize3cm\tiny\raggedright\pretolerance10000
\noindent #1\hfill}\hss}\vbox to8pt{\vfil}\vss}}}

\newcounter{ssig}
\newcounter{ttig}
\title[A generalized integrability problem for $G$-structures]
{A generalized integrability \\ problem for $G$-structures}
\author{Andrea Santi}
\address{
Andrea Santi,
Dipartimento di Matematica e Informatica, Universit\'a di Parma, Parco Area delle Scienze
$53$/A, $43124$, Parma, Italy}
\email{asanti.math@gmail.com, andrea.santi@unipr.it}
\thanks{{\it Acknowledgments}. This research was partially supported by the Project 
Firb 2012 {\it Geometria differenziale e teoria geometrica delle funzioni}, by GNSAGA
of INdAM and by project F1R-MTH-PUL-08HALO-HALOS08 of University of Luxembourg.}
\begin{document}
\newcommand{\bnabla}{\mbox{\boldmath$\nabla$}}
\begin{abstract} 
Given an $\wt n$-dimensional manifold $\wt M$ equipped with a $\wt G$-structure $\wt\pi:\wt P\rightarrow \wt M$, there is a naturally induced $G$-structure $\pi: P\rightarrow M$ on any submanifold $M\subset\wt M$ that satisfies appropriate regularity conditions. We study generalized integrability problems for a given $G$-structure $\pi: P\rightarrow M$, namely the questions of whether it is locally equivalent to induced $G$-structures on regular submanifolds of homogeneous $\wt G$-structures $\wt\pi:\wt P\to \wt{H}/\wt{K}$.  
If $\wt\pi:\wt P\to \wt{H}/\wt{K}$ is flat $k$-reductive
we introduce a sequence of generalized curvatures taking values in appropriate cohomology groups and prove that the vanishing of these curvatures are necessary and sufficient conditions for the solution of the corresponding generalized integrability problems.
\end{abstract}
\maketitle
\null \vspace*{-.25in}
\section{Introduction}
\label{Introduction}
\setcounter{section}{1}
\setcounter{equation}{0}
\vskip0.3cm
Let $M$ be an $n$-dimensional manifold and $\pi:L(M)\rightarrow M$ the linear frame bundle of $M$, that is 
the $\GL_{n}(\bR)$-bundle of all linear frames $u=(e_i)$ of the tangent spaces of $M$. We recall that, given a closed subgroup $G\subset\GL_{n}(\bR)$, a {\it $G$-structure on $M$} is a $G$-reduction $P$ of  $L(M)$, a bundle formed by a collection of linear frames with the property that any two of them $u=(e_i)$, $u'=(e_i')$ at a point $x\in M$ 
satisfy an equality of the form
$$
e_j'=\sum_{i=1}^{n} A^i_je_i\,,\qquad\text{for some}\quad A=(A_j^i)\in G\ 
$$
(see, e.g., \cite{Kb, St}).
Any $G$-structure is naturally endowed with the so-called {\it soldering form} $\vartheta:TP\rightarrow \bR^n$, defined by
$
\vartheta_u(v)=(v^1,\dots,v^n)
$,
where the $v^i$ are the components of the vector $\pi_*(v)\in T_{\pi(u)}M$ w.r.t. the frame $u=(e_i)$.

One of the simplest examples of a $G$-structure is given by the flat $G$-structure $P^o$ on $\bR^n$, the collection of all linear frames $(e_i)$ of the form
$$
e_i=\sum_{j=1}^n A_i^j\left.\frac{\partial}{\partial x^j}\right|_x\;,\;\;\;(A^j_i)\in G\;,\;\;\;(x^i)=\text{standard coord. of}\;\bR^n\,.
$$

An arbitrary $G$-structure $\pi:P\rightarrow M$ is called {\it (locally) integrable} if it satisfies the following: around any $x\in M$ there exists a local diffeomorphism $f:\cU\rightarrow \bR^n$ such that the bundle $f_*(P|_{\cU})$ of pushed-forward frames coincides with $P^{o}|_{f(\cU)}$.
Equivalently one may say that $P$ is locally integrable if $M$ can be covered by charts $\xi=(y^i)$ with the property that the coordinate frames $(\frac{\partial}{\partial y^i}\big|_x)$ are in $P$ at all points.
The problem of determining whether a given $G$-structure $\pi:P\rightarrow M$ is integrable or not is usually called the {\it integrability problem}. When $G\subset\GL_n(\bR)$ is of {\it finite type}, such problem is completely solvable. Indeed it turns out that, in this case, $P$ is integrable if and only if all elements of a special collection of $G$-equivariant maps in the
Spencer cohomology groups $H^{p,2}(\gg)$, $\gg=Lie(G)$, are trivial (see \cite{G}).
\medskip\par
Another natural construction of $G$-structures is the following. Let $\wt\pi:\wt P\rightarrow \wt M$ be a $\wt G$-structure on a manifold $\wt M$ of dimension $\wt n=n+m$ and $M\subset \wt M$ a submanifold of dimension $n$. We say that $M$ is {\it $\wt P$-regular} if the collection
of frames $u\in\wt P$ that are {\it adapted} to $M$ 
(i.e. with $x=\wt\pi(u)\in M$ and the first $n$ vectors tangent to $M$) 
constitutes a principal bundle $\pi_\sharp=\wt\pi|_{P_\sharp}:P_\sharp\rightarrow M$ with structure group
\beq
\label{semi1}
G_{\sharp}=\left\{\wt A=\begin{pmatrix} A & B \\ 0 & D\end{pmatrix}\in\wt G\;|\;A\;\text{size}\;n\times n\right\}.
\eeq
Note that the quotient $G_\sharp/N_\sharp$ of $G_\sharp$ by the closed and normal subgroup 
\beq
\label{semi2}
N_{\sharp}=\left\{\wt A=\begin{pmatrix} I & B \\ 0 & D\end{pmatrix}\in\wt G\;\right\}
\eeq 
is isomorphic to
$
G=\left\{A\in \GL_{n}(\bR)\,|\,A=\text{upper left block of some}\,\wt A\in G_\sharp \right\}
$.
If $M\subset \wt M$ is $\wt P$-regular, the quotient bundle $\pi:P=P_\sharp/N_\sharp\rightarrow M$ with structure group $G=G_\sharp/N_\sharp$ is naturally identifiable with the $G$-structure on $M$ given by the frames formed by the first $n$ vectors of the adapted frames.
We call it $G$-{\it structure induced by $\wt P$}; its
frames $u\in P$ are the {\it induced frames}.
First examples of such induced structures are the orthonormal frame bundles of $n$-dimensional submanifolds of the flat Riemannian manifold $\bR^{n+m}$. 

Now, given a $G$-structure $\pi:P\rightarrow M$ and an homogeneous $\wt G$-structure $\wt\pi:\wt P\rightarrow \wt M$ on a manifold of dimension $\wt n=n+m$, one may ask whether around any $x\in M$ there is a local embedding $f:\cU\subset M\rightarrow \wt M$ 
such that $f(\cU)$ is $\wt P$-regular and the pushed-forward bundle $f_*(P|_\cU)$ is the $G$-structure induced by $\wt P$. 
If this occurs, we say that $f$ is a {\it $(P,\wt P)$-regular embedding}
and that $P$ is {\it locally immersible into $(\wt M, \wt P)$} (shortly, {\it immersible in $\wt P$}). By analogy with the classical integrability problem, we call the question on existence of such local immersions the {\it generalized integrability problem}.
\medskip
\par
A crucial tool to deal with a given generalized integrability problem is provided by the following theorem. 

Let $\wt M=\wt{H}/\wt{K}$ be an $\wt n$-dimensional homogeneous manifold with effective action and $\wt\pi:\wt P\rightarrow \wt M$ one of its associated homogeneous $\wt G$-structures. Each of them is a subbundle $\wt P\subset L(\wt M)$ given by all frames of an orbit in $L(\wt M)$ of the group of diffeomorphisms $\wt H\subset \operatorname{Diff}(\wt M)$ (see \S \ref{gradisca} for details). Let also $\tgh=Lie(\wt H)$, $\tgk=Lie(\wt K)$ and denote by $\wt p:\wt{\gh}\rightarrow\wt{\gh}/\wt{\gk}\simeq T_o\wt M$, $o=e\wt K$, the canonical projection onto $\wt{\gh}/\wt{\gk}$. In the following, the space $T_o\wt M$ is identified with $\bR^{\wt n}$ by means of a linear frame $u_o\in\wt P|_o$.
\begin{theorem}
\label{topolino}
Let $\pi:P\rightarrow M$ be a $G$-structure on an $n$-dimensional manifold $M$ with soldering form $\vartheta=(\vartheta^1,\ldots,\vartheta^n)$. 
Then $P$ is locally immersible in $(\wt M=\wt H/\wt K,\wt P)$ if and only if 
\begin{itemize}
\item[a)]  $G\simeq G_{\sharp}/N_\sharp$, where $G_\sharp, N_\sharp\subset\wt G$ are as in \eqref{semi1} and \eqref{semi2} respectively;
\item[b)] for any local section $s:\cU\subset M\rightarrow P$, there exists an $\tgh$-valued local $1$-form $\omega$ on $M$ which satisfies 
\beq 
\label{simbad1}
\wt p\circ\omega=(s^*\vartheta^1,\ldots,s^*\vartheta^n,0,\ldots,0)
\eeq 
and
\beq
\label{simbad2}
d\omega+\frac{1}{2}[\omega,\omega]=0\ .
\eeq
\end{itemize}
Furthermore, if $(b)$ holds for a given section $s:\cU\rightarrow P$, then it holds for any other section on the same open set $\cU$.
\end{theorem}
The proof of Theorem \ref{topolino} is simple. The necessity is a consequence of the fact that, if there is the required immersion $f:\cU\subset M\rightarrow\wt M$, one can pull-back the Maurer-Cartan form $\omc$ of $\wt H$ on $M$ and obtain a $1$-form $\omega$ that satisfies \eqref{simbad1} and \eqref{simbad2}. The sufficiency is obtained by an appropriate use of Frobenius' Theorem (see \S \ref{firstorder}). 
\par\medskip
Theorem \ref{topolino} gives an unified approach for determining complete sets of necessary and sufficient conditions for the existence of local immersions into homogeneous spaces. To have a flavour of the meaning of equations \eqref{simbad1} and \eqref{simbad2}, consider the case  where 
there is a reductive decomposition
$
\tgh=\wt\gk+\wt\gm
$
and two decompositions $\tgk=\gg+\gg^\perp$ and $\wt\gm=\gm+\gm^\perp$ into the direct sum of modules for the adjoint action of $\gg=Lie(G)$. 
In this case the $\wt\gh$-valued $1$-form $\o$ has the form 
$$\o=\o^{\gg}+\o^{\gg^\perp}+\o^{\gm}+\o^{\gm^\perp}$$
and condition  \eqref{simbad1} means that $\o^{\gm}$ coincides with $s^*\vartheta$ while $\o^{\gm^\perp}$ vanishes. At the same time, condition \eqref{simbad2} corresponds to the existence of a $1$-form $\o^{\tgk}=\o^{\gg}+\o^{\gg^\perp}$ satisfying the system of equations
\beq\notag
\begin{split}
[\o^{\gg},\o^{\gm}]+[\o^{\gg^\perp},\o^{\gm}]=-d\o^{\gm}-\frac{1}{2}[\o^{\gm},\o^{\gm}]_{\wt{\gm}}\,\,\,,\\
d\o^{\gg}+d\o^{\gg^\perp}+\frac{1}{2}[\o^{\gg},\o^{\gg}]+[\o^{\gg},\o^{\gg^\perp}]+\frac{1}{2}[\o^{\gg^\perp},\o^{\gg^\perp}]=-\frac{1}{2}[\o^{\gm},\o^{\gm}]_{\wt{\gk}}\ .
\end{split}
\eeq
In the case of isometric immersions into spaces of constant curvature, these equations turn out to be equivalent to the classical Gauss-Codazzi-Ricci equations
and the property that they are the only obstructions to be satisfied is a direct consequence of Theorem \ref{topolino} and of its main corollaries (see \S \ref{Rie}).
\par\smallskip
In general, equation \eqref{simbad2} is quite involved. 
But, if $\textstyle\tgh=\sum_{p=-1}^{k-1}\wt\gh^p$ is a {\it (quasi-)graded Lie algebra of depth $1$}, it splits into a sequence of equations that can be solved with an iterative process (see \S\ref{kreductive} and \S\ref{simdor}). Further, if $\wt P$ is integrable and of finite type, in analogy with the results in \cite{G} there is a sequence of locally defined maps on $\cU\subset M$ which take values in appropriate cohomology groups and vanish if and only if condition (b) of Theorem \ref{topolino} is satisfied, i.e. if and only if the generalized integrability problem is solvable (see Theorem \ref{r-1}). 
These maps are called {\it essential $(G,\wt G)$-curvatures} and they generalize, for instance, the classical Riemannian curvature in the integrability problem of Riemannian metrics. Under appropriate cohomological conditions, the essential $(G,\wt G)$-curvatures are maps which are  globally well-defined on the total space of $\pi:P\to M$ and $G_\sharp$-equivariant (see Proposition \ref{finesez} and Theorem \ref{quaternion}). 
\medskip\par
Several applications of the general theory of $\wt G$-structures $\wt\pi:\wt P\to \wt M$ and the $G$-structures $\pi:P\to M$ induced on $\wt P$-regular submanifolds $M\subset \wt M$ are considered in \S\ref{teoria}. There we first obtain in a unified fashion the classical Gauss-Codazzi-Ricci equations for the isometric immersions of a Riemannian manifold $(M,g)$ into Euclidean spaces, spheres and hyperbolic spaces. Then we consider the case of local conformal immersions of $(M,g)$ into spheres and get a new proof for Akivis' theorem on existence of such immersions (\cite{AK61, AG}). Additional new results on this topic are also given.
Finally we show how the general results yield the classical result by Andreotti and Hill on the existence of CR local embeddings of CR manifolds in $\bC^{\wt n}$ (\cite{AH}) and we discuss the theorem on the existence of local CR quaternionic embeddings of an almost CR quaternionic manifold in $\bH$P$^{\wt n}$
that we prove in \cite{ASQ}.
\medskip
\par
The paper is organized as follows. In \S 2 we give the basic definitions and properties of homogeneous $\wt G$-structures and induced $G$-structures and prove Theorem \ref{topolino}. In \S 3 we consider the case of (quasi-)graded Lie algebras, describe the above mentioned iterative process and prove Theorem \ref{r-1}, the main result of this section. In \S\ref{tortaele} and \S 5 we restrict to local immersions in flat $2$-reductive $\wt G$-structures and present the previously mentioned examples of applications.
\medskip\par
Before concluding, we have to mention that a general study of the generalized integrability problem was also done by Rosly and Schwarz in \cite{RS}. In fact, this paper was strongly inspired by that one. There, the authors considered immersions into integrable $\wt G$-structures of finite type, estabilished Theorem \ref{topolino} and determined a finite
sequence of obstructions. Main novelties of this paper with respect to \cite{RS} are the proof of Theorem \ref{topolino} in coordinate-free language, which allows immediate 
extensions to the case of local immersions into homogeneous manifolds, and the introduction of a new sequence of cohomology groups and the corresponding essential curvatures. Note that this new setting can also be easily adapted to deal with other kind of $G$-structures, as for instance those of higher order (see \cite{ASP}). We also remark that the essential curvatures considered in this paper are in general different from Rosly and Schwarz' curvatures and actually determine a set of obstructions which completes the one presented  in \cite{RS} (see Remark \ref{cisiamo!!}).
\vskip0.2cm\par\noindent
{\it Acknowledgements.} 
Part of this work was done while the author was a post-doc at the University of Parma. 
The author would like to thank the Mathematics Department and in particular A. Tomassini and 
C. Medori for support and ideal working conditions.
The author is also grateful to A. Spiro for useful discussions on various aspects of this paper
and to the anonymous referee for his helpful comments and suggestions.
\par\medskip
\section{$G$-structures immersible into homogeneous $\wt G$-structures}
\label{firstorder}
\setcounter{section}{2}
\setcounter{equation}{0}
\par\noindent
\subsection{\it Homogeneous $G$-structures.}\label{gradisca}
Let $M=H/K$ be an $n$-dimensional homogeneous manifold with an effective action. We throughout assume that $H$ is connected and simply connected and $K$ connected. For any $h\in H$, we denote
by $L_{h}\in\operatorname{Diff}(M)$ the corresponding left action.

Given a frame $u_o=(e_i)\in L(M)$ at $o=eK$, let 
$$P=H\cdot u_o=\{h\cdot u_o:=(L_h{}_* e_i)\,,\,h\in\ H\}\subset L(M)$$ 
be the $H$-orbit of $u_o$. One gets in this way a reduction $\pi:P\rightarrow M$ of $L(M)$, which we call the {\it homogeneous $G$-structure associated with $M=H/K$}. 
\smallskip
\par
Note that if $K_1$ is the kernel of the isotropy representation
$$
K\rightarrow \GL(T_oM)\;,\qquad k\mapsto L_{k}{}_*|_{o}\,,
$$
then the map $h\mapsto h\cdot u_o$
determines a natural bundle isomorphism  between $H/K_1\rightarrow M=H/K$ and the $G$-structure $\pi:P\rightarrow M$, with $G\simeq K/K_1$. 

From now on, we assume a fixed choice of $u_o$ and tacitly use it to identify $T_oM$ with $\bR^n$ and $\pi:P\rightarrow M$ with $H/K_1\rightarrow H/K$.
Then, if we consider the canonical projections
$$\pi_1:H\rightarrow P=H/K_1\;\;,\;\;\;\;p:\gh\rightarrow \gh/\gk\simeq T_oM\,,$$
the soldering form $\vartheta$ of $P$ can be determined from the Maurer-Cartan form $\varpi$ of $H$ as the unique $\bR^n$-valued $1$-form which satisfies
\beq
\label{mizz}
\pi_1^*\vartheta=p\circ\varpi\ .
\eeq
\vskip0.2cm\par
In the following, we define {\it homogeneous $G$-structure} any pair $(H/K,P)$ consisting of an homogeneous manifold $M=H/K$ and a fixed associated homogeneous $G$-structure $P$. Further, an homogeneous $G$-structure $(H'/K',P')$ is defined {\it extension of $(H/K,P)$} if $H'/K'$ is locally diffeomorphic to $H/K$, $P'$ locally equivalent to $P$ as a $G$-structure and $\gh=Lie(H)$ is a proper subalgebra of $\gh'=Lie(H')$ with $\gk=\gh\cap \gk'$. The homogeneous $G$-structure $(H/K,P)$
is {\it inextendible} if it does not admit any extension.
\smallskip\par
As mentioned in the introduction, Theorem \ref{topolino} fully characterizes immersible $G$-structures. 
Before the proof, we give a simple characterization of $\wt P$-regular submanifolds, which hold for any (not necessarily homogeneous) $\wt G$-structure $\wt\pi:\wt P\to \wt M$.
\begin{lemma}\label{characterization}
A submanifold $M\subset \wt M$ is $\wt P$-regular if and only if for any $x\in M$ there exists an adapted frame $u$ at $x$. In this case the principal bundle $\pi_\sharp:P_\sharp\to M$ is a $G_\sharp$-reduction  of $\wt P|_M$ with $G_{\sharp}$ as in \eqref{semi1}.
\end{lemma}
This lemma is proved by first checking the existence of local smooth
sections of $\wt P|_M$ taking values in $P_\sharp$ and then applying Lemma 1 in \cite[p. 84]{KN1}.
We omit the details for the sake of brevity.
\subsection{\it Proof of Theorem \ref{topolino}}
Let $f:\cU\subset M\rightarrow \wt M$ be a $(P,\wt P)$-regular embedding, which we use to identify $\cU$ with a submanifold of $\wt M$ and $P|_\cU$ with the induced $G$-structure. Without loss of generality we may assume $\cU$ to be contractible,
consider a local section $s_\sharp:\cU\rightarrow P_\sharp$ of the subbundle $P_\sharp\subset\wt P|_\cU$ and the corresponding section $s:\cU\rightarrow P$ of $P$  induced via the identification $P\simeq P_\sharp/N_\sharp$. 
By local triviality of the fibration $\tilde\pi_1:\wt{H}\rightarrow \wt P=\wt H/\wt K_1$, we may also pick  a local section $\wt s:\cU\rightarrow \wt{H}|_\cU$ satisfying $\wt\pi_1\circ \wt s=s_\sharp$. The situation is summarized by the following diagram: 
\centerline{\hskip2cm
 \begin{picture}(300,215)(0,0)
 \setlength{\unitlength}{1pt}
 \linethickness{0.1mm}
\put(130, 190){$\wt H|_{\cU}$}
\put(141, 185){\vector(0,-1){45}}
\put(130, 125){$\wt P|_\cU=\wt H/\wt K_1|_{\cU}$}
\put(141, 120){\vector(0,-1){45}}
\put(137, 61){$\cU\subset \wt M=\wt H/\wt K$}
\put(11, 125){$P_\sharp$}
\put(23, 122){\vector(2,-1){41}}
\put(60, 90){$P=P_\sharp/N_\sharp$}
\put(94, 84){\vector(2,-1){38}}
\put(25,130){\vector(1,0){102}}
\qbezier(132,61)(20,50)(14,119)
\put(14,123){\vector(0,1){0}}
\qbezier(132,61)(95,60)(79,85)
\qbezier(145,70)(280,117)(146,183)
\put(143,186){\vector(-1,1){0}}
\put(77,89){\vector(-1,2){0}}
\put(125,160){$\wt\pi_1$}
\put(127,90){$\wt\pi$}
\put(215,102){$\wt s$}
\put(50,55){$s_\sharp$}
\put(115,77){$\pi$}
\put(78,68){$s$}
 \linethickness{0.05mm}
\put(30,134){\oval(15,8)[bl]}\put(27.5,134){\oval(10,8)[tl]}
 \linethickness{0.1mm}
\end{picture}
 }
\vskip-1.5cm\par\noindent
Let $\omega=\wt s^*\varpi$, where $\varpi$ is the Maurer-Cartan form  of $\wt H$. Then equation \eqref{simbad1} is checked observing that
$$
\wt p\circ\omega=\wt s^*(\wt p\circ\varpi)=\wt s^*(\wt\pi_1^*\wt\vartheta)=s_\sharp^*\wt\vartheta=(s^*\vartheta^1,\ldots,s^*\vartheta^n,0,\ldots,0)\,,
$$
where the second equality follows from \eqref{mizz} applied to the soldering form $\wt \vartheta$ of $\wt P$ and the last equality is a direct consequence of definitions of soldering $1$-forms and induced $G$-structures.
The equation \eqref{simbad2} follows from the Maurer-Cartan equation of $\varpi$. 
\medskip
\par
Conversely, let $\omega$ be a $1$-form satisfying \eqref{simbad1}, \eqref{simbad2} and consider the $1$-form $\omega-\varpi$ on the product $\cU\times \wt{H}$ with the associated distribution $\cD\subset T(\cU\times\wt H)$ of rank $n$, defined by
$$\cD=\cup_{y\in\cU\times\wt H}\cD_y\;\;,\;\;\;\; \mathcal{D}_{y}=\operatorname{Ker}(\omega-\varpi)|_{y}\,.$$ 
Using \eqref{simbad2}, one checks that $\cD$ is involutive, hence integrable by Frobenius' Theorem. For a chosen $x\in\cU$, an integral leaf through the point $(x,e)$ is the graph of a unique local map $\wt s:\cU\rightarrow \wt{H}$ with $\omega=\wt s^*\varpi$ and $\wt s(x)=e$.
\smallskip\par
We now claim that, by possibly taking a smaller $\cU$, the map
$$
f=\wt\pi\circ\wt\pi_1\circ \wt s:\cU\subset M\rightarrow \wt M
$$
is a $(P,\wt P)$-regular embedding.
Using \eqref{simbad1}, one sees that $f_*|_{x}$ is injective so that $f$ is a local embedding around $x$.
Moreover, given the section $s_\sharp:=\wt\pi_1\circ \wt s:\cU\rightarrow \wt P|_\cU$ of $\wt P|_\cU$,
we may consider the $G_\sharp$-reduction 
$$P_\sharp=\{s_\sharp(x)\cdot\wt A\,:\, x\in \cU,\, \wt A\in G_\sharp\}\subset\wt P|_\cU\ .$$
By \eqref{mizz} and \eqref{simbad1}, we get
$$
\hskip-0.27cm s_\sharp^*\wt\vartheta=
\wt s^*\wt\pi_1^*\wt\vartheta=
\wt p\circ (\wt s^*\varpi)=\wt p\circ\omega=(s^*\vartheta^1,\ldots,s^*\vartheta^n,0,\ldots,0)\,. 
$$
This relation between the soldering forms $\wt\vartheta$, $\vartheta$ of $\wt P$ and  $P$, respectively, says that $f$ is $(P,\wt P)$-regular
with bundle of adapted frames isomorphic with $P_\sharp$.
\smallskip\par
Since we just proved that the existence of solutions to \eqref{simbad1} and \eqref{simbad2} for a given section $s:\cU\rightarrow P$
is equivalent to  local immersibility of $P$ (a property clearly independent on the choice of $s$), the last claim of the theorem follows immediately.
\qed
\begin{remark}
Let $(\wt H'/\wt K',\wt P')$ be a proper extension of the homogeneous $\wt G$-structure $(\wt H/\wt K,\wt P)$. Since any $\wt\gh$-valued $1$-form is an $\wt\gh'$-valued $1$-form satisfying additional conditions, the equations \eqref{simbad1} and \eqref{simbad2} on $\wt\gh$-valued $1$-forms are equivalent to the system on $\wt\gh'$-valued $1$-forms given by \eqref{simbad1}, \eqref{simbad2} and such constraints. On the other hand, since $\wt P'\simeq \wt P$, Theorem \ref{topolino} implies that there is an $\wt\gh$-valued solution of \eqref{simbad1} and \eqref{simbad2} if and only if there is an $\wt\gh'$-valued solution. 
Due to this, inextendible homogeneous $\wt G$-structures are more appropriate than the extendible ones
for determining  minimal sets of obstructions for a given  generalized integrability problem.
\end{remark}
\par\medskip
\section{$k$-reductive $\wt G$-structures and $(G,\wt G)$-curvatures}
\label{sez3}
\setcounter{section}{3}
\setcounter{equation}{0}
\vskip0.25cm
\subsection{\it $k$-reductive $G$-structures}\label{kreductive}
We recall that a Lie algebra $\gh$ is a {\it graded Lie algebra of depth $1$ and height $k$} 
(shortly, {\it graded Lie algebra of height $k$}) if it admits a direct sum decomposition of vector spaces 
\beq
\label{decompA}
\gh=\sum_{p=-1}^{k-1}\gh^p
\eeq
such that 
\beq
\label{decompB}
[\gh^{r},\gh^{s}]\subset \gh^{r+s}\qquad\text{for all}\;r,s,
\eeq 
where $\gh^{r+s}=\{0\}$ for $r+s\notin\{-1,0,\ldots,k-1\}$. If a Lie algebra $\gh$ admits a vector space decomposition \eqref{decompA} and satisfies \eqref{decompB} only when $r+s \geq -1$, we call it {\it quasi-graded of height $k$}. 

Note that {\it any} Lie algebra $\gg$  can be considered as a quasi-graded Lie algebra of height $k = 0$, while the quasi-graded Lie algebras of height $k=1$ coincide with the Lie algebras endowed with a reductive decomposition.
\smallskip\par
The $n$-dimensional homogeneous spaces $M = H/K$ with an effective action, quasi-graded Lie algebra $\gh = Lie(H)$ of height $k$ and isotropy subalgebra 
$$
\gk = Lie(K)=\sum_{p=0}^{k-1}\gh^p
$$
are called {\it $k$-reductive}; they have been considered  by I. L. Kantor  (\cite{Ka, AS}) as natural  generalizations of reductive homogeneous spaces.
In this case
$$
\dim \gh^{-1} = \dim (\gh/\gk) = n\qquad \text{and}
$$
$$
\text{for any} \ X \in \gh^p, \ p \geq 0\ ,\quad [X, \gh^{-1}] = 0\quad \text{implies}\quad X = 0\ ,
$$
by transitivity and effectiveness of the action of $H$ on $M$.

We call {\it $k$-reductive $G$-structure} any homogeneous $G$-structure $(H/K,P)$ associated with a $k$-reductive space $H/K$.
Under an identification $\gh/\gk\simeq\bR^{n}$ given by a given frame $u_o$ of $P|_o$, the Lie algebra $\gg\subset\ggl_{n}(\bR)$ of the structure group $G$ of $P$ is identifiable with $\gh^0$.  

When $\gh$ is a graded Lie algebra, it is known that $(H/K, P)$ is locally identifiable with the flat $G$-structure on $\bR^{n}$, i.e. it is  {\it locally integrable}.
We call it {\it flat $k$-reductive $G$-structure}. Its Lie algebra $\textstyle\gh=\sum_{p=-1}^{k-1}\gh^{p}$ is a subalgebra of the maximal prolongation of $\gh^0 \subset \ggl_{n}( \bR)$, defined as follows. 
\begin{definition}  The {\it maximal transitive prolongation} of a subalgebra $\gh^0$ of $\ggl_{n}( \bR)$  is the maximal  graded Lie algebra 
$$
\gh_{\infty}=\sum_{p=-1}^{\infty}\gh_\infty^{p}\ ,
$$
with $\gh_\infty^{-1}= \bR^{n}$, $\gh_\infty^{0}=\gh^{0}$  and satisfying the following two properties:
\begin{itemize}
\item[1)]
the adjoint action of $\gh_\infty^{0}$ on $\gh_\infty^{-1}$ is the standard action of $\gh^{0}$ on $\bR^{n}$,
\item[2)]
for any $X\in\gh_{\infty}^p$, $p\geq 0$, if $[X,\gh_\infty^{-1}] = 0$ then $X = 0$.
\end{itemize}
\end{definition}
\smallskip\par
It is known that the maximal prolongation $\gh_\infty$ is uniquely determined by the algebra $\gh^{-1}+\gh^0  \subset\bR^{n}  +\ggl_{n}(\bR)$ and it might be infinite-dimensional. When it is finite-dimensional, the Lie algebra $\gg = \gh^0$ is called {\it of finite type}. Accordingly, also the $G$-structures with Lie algebras $Lie(G) = \gg$ are called {\it of finite type} (for further details, see e.g. \cite[Ch. 1]{Kb}). 
\smallskip\par
A flat $k$-reductive  $G$-structure $(H/K, P)$, with $P$ of finite-type and 
$\gh$ equal to its maximal prolongation $\textstyle\gh = \gh_\infty = \sum^{k-1}_{p = -1}\gh_\infty^p$ is inextendible (this follows from Sternberg prolongation theory of bundles, see \cite{St}). 
Further, one can check that  the following holds (see \cite[Thms. 17-18]{Ka}):  {\it if an  homogeneous space   $H'/K'$ admits an invariant Cartan connection modelled on $H/K$ and is of maximal dimension (that is $\dim H' = \dim H$) then: 
\begin{itemize}
\item[$-$] $H'/K'$ is $k$-reductive,  
\item[$-$] $Lie(H')$ is a quasi-graded Lie algebra with underlying graded vector space isomorphic to $\textstyle\gh$,
\item[$-$] $\textstyle Lie (K') \simeq Lie(K) = \sum_{p=0}^{k-1}\gh_\infty^p$.
\end{itemize}}\par\noindent
First examples of such homogeneous spaces and their associated $k$-reductive  $G$-structures $(H'/K', P')$ are the three simply connected Riemannian manifolds of constant curvature and their orthonormal frame bundles (see \S \ref{Rie}).
\subsection{\it Immersible $G$-structures of order $p$ and associated $(G,\wt G)$-curvatures}\label{simdor}
Let $(\wt H/\wt K,\wt P)$ be a $k$-reductive $\wt G$-structure on an $\wt n$-dimensional homogeneous manifold $\wt M=\wt H/\wt K$ and $\pi:P\rightarrow M$ a $G$-structure on an $n$-dimensional manifold $M$ with soldering form $\vartheta=(\vartheta^1,\ldots,\vartheta^n)$, where $G$ is isomorphic to $G_\sharp/N_\sharp$ with $G_\sharp, N_\sharp\subset\wt G$ as in \eqref{semi1} and \eqref{semi2}. 

Combining Theorem \ref{topolino} with the decomposition $\textstyle\wt\gh=\sum_{p=-1}^{k-1}\wt\gh^p$, one gets that $P$ is locally immersible in $\wt P$ if and only if, for any local section $s:\cU\subset M\rightarrow P$ of $P$, there exist $k+1$ $1$-forms 
$$\omega^{p}:T\cU\rightarrow \wt\gh^{p}\;,\qquad -1\leq p\leq k-1\,,$$ 
with
\beq
\label{MCC}
\omega^{\mathrm{-1}}=(s^*\vartheta^1,\ldots,s^*\vartheta^n,0,\ldots,0)
\eeq
and the others satisfying the system of $k+1$ equations
\beq
d\omega^{p-1}+\frac{1}{2}\sum_{r=0}^{p-1}[\omega^r,\omega^{p-1-r}]+\frac{1}{2}[\omega^{-1},\omega^{-1}]^{\wt\gh^{p-1}}=-[\omega^{-1},\omega^p]\;,\;\; 0\leq p\leq k-1\,,
$$
$$
\label{gradiamo}
d\omega^{k-1}+\frac{1}{2}\sum_{r=0}^{k-1}[\omega^r,\omega^{k-1-r}]+\frac{1}{2}[\omega^{-1},\omega^{-1}]^{\wt\gh^{k-1}}=0\ ,
\eeq
where we denoted by $(\cdot)^{\wt \gh^{r}}$ the natural projection onto the subspace $\wt \gh^r$ of $\wt \gh$.
\smallskip\par
We stress the fact that, if \eqref{MCC} and \eqref{gradiamo} hold for 
a choice of $s$, then they hold for any other section of $P$ on the same open set $\cU$. The very same property is also shared by the subsystem of equations which consists of \eqref{MCC} and just the first $p$ equations of \eqref{gradiamo} for some $p<k+1$. 
Due to this, in order to simplify the quest for solutions, it is always  advisable to
 make a choice for the section $s$ (from which  the associated 1-form $\o^{-1}$ is always uniquely determined by \eqref{MCC})  that makes  the explicit expression of \eqref{MCC} as simple as  possible (see for instance, \S\ref{conforme}).
\bigskip\par
All this motivates the following:
\begin{definition}
\label{acqua}
The $G$-structure $P$ is defined {\it locally immersibile into $(\wt H/\wt K,\wt P)$ to order $p$} if, around any point, it admits a local section $s$ and forms $\omega^{-1},\dots,\omega^{p-1}$ which satisfy \eqref{MCC} and the first $p$ equations of \eqref{gradiamo}. Any $p+1$-tuple
$(s, \o^0, \ldots, \o^{p-1})$ as above is called {\it admissible}.  
\end{definition}
Let $P$ be immersible to order $p$ and $(s, \o^0, \ldots, \o^{p-1})$ an admissible tuple.
Let also $V=\bR^{\wt n}\simeq \wt\gh^{-1}$ and denote by $W\simeq\bR^n$ the $n$-dimensional subspace of $V$ which corresponds to the vanishing of the last $\wt n -n$ standard coordinates of $\bR^{\wt n}$. 
For any integer $r=0,1,\ldots,p$ we define {\it annihilator  of  level $r$ in $\wt\gh^{p-1}$} 
the subspace of $\wt\gh^{p-1}$ given by
$$
\gc^{p-1}_{r}=\left\{X\in \wt\gh^{p-1}\,|\,[\ldots[[X,w_1],w_2]\ldots,w_{r}]=0\;\;\text{for all}\;\;w_1,\ldots,w_r\in W\right\}\ .
$$
For instance $\gc^{p-1}_0=0$ and $\gc^{p-1}_1$ is the centralizer of $W$ in $\wt\gh^{p-1}$. 
We also set $\gc^{p-1}_{p+1}:=\wt\gh^{p-1}$ and consider the increasing filtration of $\wt\gh^{p-1}$ given by
$$
0=\gc^{p-1}_0\subset\gc^{p-1}_1\subset\cdots\subset\gc^{p-1}_r\subset\cdots\subset \gc^{p-1}_{p}\subset\gc^{p-1}_{p+1}=\wt\gh^{p-1}\ .
$$
Let $p'_{r}: \wt\gh^{p-1}\to \wt\gh^{p-1}/\gc^{p-1}_r$ be the natural projection.
We define {\it  $(G, \wt G)$-curvature of level $r$ and order $p+1$ associated with $(s, \o^0, \ldots, \o^{p-1})$} the $\wt\gh^{p-1}/\gc^{p-1}_r$-valued $2$-form on the domain of $s$ given by
\beq
\label{tavolo}
\wh\Omega^{p-1}_r=p'_r\circ (d\omega^{p-1}+\frac{1}{2}\sum_{s=0}^{p-1}[\omega^s,\omega^{p-1-s}]+\frac{1}{2}[\omega^{-1},\omega^{-1}]^{\wt\gh^{p-1}})\ .
\eeq

For any $r=0,1,\ldots,p$ we fix a subspace $\gc^{p-1\perp}_r$ complementary to $\gc^{p-1}_r$ in $\gc^{p-1}_{r+1}$. One gets in this way a vector space direct sum decomposition
\beq
\label{complementare}
\wt\gh^{p-1}=\bigoplus_{s=0}^{p}\gc^{p-1\perp}_s\;\;,\;\;\;\text{where}\;\;\;\gc^{p-1}_r=\bigoplus_{s=0}^{r-1}\gc^{p-1\perp}_s\;\;,
\eeq
which determines a corresponding projection $p''_r: \wt\gh^{p-1} \to \gc^{p-1\perp}_r$. 
We define {\it  complementary $(G, \wt G)$-curvature of level $r$ and order $p+1$ associated with $(s, \o^0, \ldots, \o^{p-1})$ and the decomposition \eqref{complementare}} the $\gc^{p-1\perp}_r$-valued $2$-form on the domain of $s$ given by
\beq
\label{tavolino}
\wt\Omega^{p-1}_r=p''_r\circ(d\omega^{p-1}+\frac{1}{2}\sum_{s=0}^{p-1}[\omega^s,\omega^{p-1-s}]+\frac{1}{2}[\omega^{-1},\omega^{-1}]^{\wt\gh^{p-1}})\ .
\eeq

Finally, we define {\it total $(G, \wt G)$-curvature of order $p+1$ associated with $(s, \o^0, \ldots, \o^{p-1})$} the $\wt\gh^{p-1}$-valued $2$-form on the domain of $s$ given by
$$
\Omega^{p-1}=d\omega^{p-1}+\frac{1}{2}\sum_{r=0}^{p-1}[\omega^r,\omega^{p-1-r}]+\frac{1}{2}[\omega^{-1},\omega^{-1}]^{\wt\gh^{p-1}}\;.
$$
We observe that for any fixed choice of a decomposition \eqref{complementare} as above,
$\wt\gh^{p-1}/\gc^{p-1}_r$ is naturally identified with $\displaystyle\bigoplus_{s=r}^p\gc^{p-1\perp}_s$ and the total curvature is the sum 
$$
\O^{p-1}=\wh\O^{p-1}_r+\sum_{s=0}^{r-1}\wt\O^{p-1}_{s}
$$
of the curvature of level $r$ and the complementary curvatures of lower levels.
\smallskip\par
The question whether a given $G$-structure $P$ is locally immersible in $\wt P$ is therefore solvable
by the following iterative procedure. 
First, fix a section $s:\cU\subset M\to P$ and set
$\omega^{-1}=(s^*\vartheta^1,\ldots,s^*\vartheta^n,0,\ldots,0)$. Secondly, look for an $\wt\gh^0$-valued $1$-form $\omega^0$ which solves the first equation of \eqref{gradiamo}. Such equation is purely algebraic in $\omega^0$ and, using the total $(G,\wt G)$-curvature $\Omega^{-1}$, can be written as 
$$
[\omega^{-1},\omega^{0}]=-\Omega^{-1}\ .
$$
If there exists a solution $\omega^0$, $P$ is immersible to order $1$. We may then compute the
total $(G,\wt G)$-curvature $\Omega^0$ and look for a $\wt\gh^1$-valued $1$-form $\omega^1$ solution to the next (algebraic) equation 
$$
[\omega^{-1},\omega^{1}]=-\Omega^{0}\ ,\;\;
$$
and so on. 
In the final step, we have that
$P$ is immersible if and only if there exists a choice of an admissible $k+1$-tuple $(s,\omega^0,\ldots,\omega^{k-1})$ whose associated total $(G,\wt G)$-curvature $\Omega^{k-1}$ vanishes identically. 
\par 
These operations shows that, at the $p^{\text{th}}$-step, the immersibility to the higher order depends on an algebraic obstruction, corresponding to the pointwise solvability of the non-homogeneous linear equation 
$$
[\omega^{-1},\omega^p]=-\Omega^{p-1}
$$
in the unknown $\omega^p$. Moreover, using the $(G,\wt G)$-curvatures of level r and the complementary curvatures, such equation 
splits for any $r=0,1,\ldots,p$ into the system of $r$ equations
\begin{align*}
&p'_r\circ[\omega^{-1},\omega^{p}]=-\wh\Omega^{p-1}_r\;\;\;\;\text{and}\\
&p''_s\circ[\omega^{-1},\omega^{p}]=-\wt\Omega^{p-1}_s\;\;\;\;\text{for any}\;s=0,\ldots,r-1\;.
\end{align*}
We respectively call them {\it equation of level $r$} and {\it complementary equations of levels $0,\ldots,r-1$}. We remark that the equation of level $r$ 
depends just on $p'_r\circ\omega^{p}$ 
and not on the entire $1$-form $\omega^p$.
\subsection{\it Immersible $G$-structures and flat $k$-reductive $\wt G$-structures}\label{biscotto}
When $(\wt H/\wt K,\wt P)$ is {\it flat $k$-reductive} the algebraic conditions on the total $(G,\wt G)$-curvature $\Omega^{p-1}$ and its components $\wh\Omega^{p-1}_r$ of level $r$ can be expressed in cohomological terms. For this, we need the following. 
\begin{lemma}
\label{GBI}
If $(\wt H/\wt K, \wt P)$ is flat $k$-reductive, the total $(G,\wt G)$-curvature associated with an admissible tuple $(s, \omega^{0},\ldots,\omega^{p-1})$ satisfies the generalized Bianchi identity
$[\omega^{-1},\Omega^{p-1}]=0$. In particular its component of level $r$ satisfies $p'_r\circ[\omega^{-1},\wh\Omega^{p-1}_r]=0$.
\end{lemma}
\begin{proof}
Since $(s, \omega^{0},\ldots,\omega^{p-1})$ is admissible and $\wt\gh^{-1}$ is abelian one has
$$d\o^{r}+\frac{1}{2}\sum_{s=0}^{r}[\omega^s,\o^{r-s}]=-[\o^{-1},\o^{r+1}]\;\;,\qquad \text{for all}\;\;-1\leq r\leq p-2\ .$$
Applying $d$ to both sides when $r=p-2$ one gets
$$
0=-d(d\omega^{p-2})=\frac{1}{2}\sum_{r=0}^{p-2}d[\omega^r,\omega^{p-2-r}]+d[\o^{-1},\o^{p-1}]\ .
$$
This equality together with admissibility of $(s, \omega^{0},\ldots,\omega^{p-1})$ yields
\begin{align*}
0&=\frac{1}{2}\sum_{r=0}^{p-2}[d\omega^r,\omega^{p-2-r}]-\frac{1}{2}\sum_{r=0}^{p-2}[\omega^r,d\omega^{p-2-r}]+[d\o^{-1},\o^{p-1}]-[\o^{-1},d\o^{p-1}]\\
&= \sum_{r=-1}^{p-2}[d\omega^{r},\o^{p-2-r}]-[\o^{-1},d\omega^{p-1}]\\
&=-\frac{1}{2}\sum_{r=0}^{p-2}\sum_{s=0}^{r}[[\omega^s,\o^{r-s}],\o^{p-2-r}]-\!\!\!\sum_{r=-1}^{p-2}[[\o^{-1},\o^{r+1}],\o^{p-2-r}]-[\o^{-1},d\omega^{p-1}].
\end{align*}
On the other hand, it is not difficult to see
$\displaystyle
\sum_{r=0}^{p-2}\sum_{s=0}^{r}[[\omega^s,\o^{r-s}],\o^{p-2-r}]=0
$
and
\begin{align*}
\sum_{r=-1}^{p-2}[[\o^{-1},\o^{r+1}],\o^{p-2-r}]&=
\frac{1}{2}\sum_{r=-1}^{p-2}[\omega^{-1},[\omega^{r+1},\o^{p-2-r}]]\\
&=\frac{1}{2}\sum_{s=0}^{p-1}[\omega^{-1},[\omega^{s},\o^{p-1-s}]]\ ,
\end{align*}
so that $\displaystyle 0=-\frac{1}{2}\sum_{s=0}^{p-1}[\omega^{-1},[\omega^{s},\o^{p-1-s}]]-[\o^{-1},d\omega^{p-1}]=-[\omega^{-1},\Omega^{p-1}]$.

It remains to prove the second part of the lemma. From $[\omega^{-1},\Omega^{p-1}]=0$ one first infers
$p'_r\circ[\omega^{-1},\Omega^{p-1}]= 0$. Since $\ad_w(\gc^{p-1}_r)\subset\gc^{p-2}_{r-1}\subset \gc^{p-2}_r$ for all $w\in W$, the claim follows directly from the fact that $\omega^{-1}$ is $W$-valued.  
\end{proof}
\smallskip\par
This lemma suggests to consider the following differential complex. For any non-negative integers $r$, $p$ and $q$ we set
$$C^{0,q}_r(\wt\gh,W):=V\otimes\Lambda^q W^{*}\;\;,\;\;\;\;   C^{p+1,q}_r(\wt\gh,W):=(\wt\gh^p/\gc_r^p)\otimes\Lambda^q W^{*}\ .$$
Consider the generalized Spencer operator $\partial:C^{p,q}_r(\wt\gh,W)\rightarrow C^{p-1,q+1}_r(\wt\gh,W)$ given by
$$
\partial c(w_1,\dots,w_{q+1}):=\sum_{i=1}^{q+1}(-1)^{i}p'_r\circ[c(w_1,\dots,w_{i-1},\hat{w_i},w_{i+1},\dots,w_{q+1}),w_i]
$$
for any $c\in C^{p,q}_r(\wt\gh,W)$ and $w_1,\dots,w_{q+1}\in W$. One can directly check that $\partial^2=0$ and consider the corresponding differential complex
$$
\cdots\overset{\partial}{\longrightarrow} C^{p+1,q-1}_r(\wt\gh,W)\overset{\partial}{\longrightarrow} C^{p,q}_r(\wt\gh,W)\overset{\partial}{\longrightarrow} C^{p-1,q+1}_r(\wt\gh,W)\overset{\partial}{\longrightarrow}\cdots\;.
$$
Let $Z^{p,q}_r(\wt\gh,W)$ and $B^{p,q}_r(\wt\gh,W)$ be the spaces of $(p,q)$-cocycles and $(p,q)$-coboundaries of this complex and set 
$$H^{p,q}_r(\wt\gh, W)=Z^{p,q}_r(\wt\gh,W)/B^{p,q}_r(\wt\gh,W)\;.$$
We call it {\it $(p,q)$-cohomology group of level $r$}. 
\vskip0.2cm\par
The following theorem extends Guillemin's results on the integrability problem for $G$-structures (see \cite{G}) to the generalized integrability problem.
\begin{theorem}
\label{r-1}
Let $(\wt H/\wt K, \wt P)$ be a flat $k$-reductive $\wt G$-structure on an homogeneous manifold $\wt M=\wt H/\wt K$ of dimension $\wt n$ and $\pi:P\rightarrow M$ a $G$-structure on a manifold $M$ of dimension $n$ with $G$ isomorphic to $G_\sharp/N_\sharp$, where $G_\sharp, N_\sharp\subset\wt G$ are as in \eqref{semi1} and \eqref{semi2}. 

Assume also that $P$ is locally immersible in $\wt P$ to order $p$ and, for any admissible tuple $(s,\omega^0,\ldots,\omega^{p-1})$ on some open set $\cU\subset M$, use the frames $s_x$, $x\in\cU$, to identify the tangent spaces $T_xM$ with $W$. Then 
\vskip0.35cm\par\noindent
i) for any admissible $(s,\omega^0,\ldots,\omega^{p-1})$, 
the associated total $(G,\wt G)$-curvature tensors
$\Omega^{p-1}|_x\in \wt\gh^{p-1}\otimes\Lambda^{2}T^*_xM\simeq C^{p,2}(\wt\gh,W)$
satisfy 
\beq
\label{pioggia}
\partial (\Omega^{p-1}|_x)=0\ 
\eeq
and  $\partial(\wh\Omega^{p-1}_r|_x)=0$ for all projections $\wh\Omega^{p-1}_r|_x\in C^{p,2}_r(\wt\gh,W)$, $1\leq r\leq p$;
\vskip0.35cm\par\noindent
ii) $P$ is locally immersible in $\wt P$ to order $p+1$ if and only if, around any point $x_o$, there exists an admissible tuple $(s,\omega^0,\ldots,\omega^{p-1})$ such that $[\Omega^{p-1}|_x]\in H^{p,2}(\wt\gh,W)$ is zero at any $x$ where it is defined.   
\end{theorem}
\begin{proof}
(i) We first observe that given an admissible tuple $(s,\omega^0,\ldots,\omega^{p-1})$ on $\cU\subset M$,
any local $\wt\gh^{s}$-valued $q$-form is naturally identified with a map from $\cU$ to $C^{s-1,q}(\wt\gh,W)$. The equation \eqref{pioggia} is then a consequence of definitions, Lemma \ref{GBI} and the fact that 
$\omega^{-1}=(s^*\vartheta^1,\ldots,s^*\vartheta^n,0,\ldots,0)$ while the last claim follows from the fact that the natural projection from $C^{p,q}(\wt\gh,W)$ to $C^{p,q}_r(\wt\gh,W)$ is a morphism of differential complexes.
\smallskip\par\noindent
(ii) We recall that $P$ is locally immersible in $\wt P$ to order $p+1$ if and only if, around any point $x_o$, there is an admissible tuple $(s,\omega^0,\ldots,\omega^{p-1})$ with associated total $(G,\wt G)$-curvature $\Omega^{p-1}$ and a smooth map
$\omega^p:\cU\rightarrow C^{p+1,1}(\wt\gh,W)$ such that $\partial\omega^p=-\Omega^{p-1}$.
%
If this is the case, the class $[\Omega^{p-1}|_x]\in H^{p,2}(\wt\gh,W)$ is zero at any point. Conversely if $[\Omega^{p-1}|_x]$ is constantly zero then, at any fixed point $x$, there is an $\omega^p_x\in C^{p+1,1}(\wt\gh,W)$ with $\partial\omega^p_x=-\Omega^{p-1}|_x$. Standard arguments yield that there exists a smooth $\wt\gh^p$-valued $1$-form $\omega^p$ on $\cU$ such that $\partial\omega^p=-\Omega^{p-1}$.
\end{proof}
In view of this result we define {\it essential $(G,\wt G)$-curvature of order $p+1$} the map
$$
\cR^{p+1}:\cU\rightarrow H^{p,2}(\wt\gh,W)\;,\;\;\;\cR^{p+1}_x:=[\Omega^{p-1}|_x]\;.
$$
It may be considered as a generalization of the $(p+1)^{\text{th}}$-order structure function introduced in \cite{G} for the classical integrability problem. 

The following corollary is a direct consequence of Theorem \ref{r-1}.
\par\noindent
\begin{corollary}\label{aereosol}
Let $(\wt H/\wt K,\wt P)$ be a flat $k$-reductive $\wt G$-structure for the graded Lie algebra $\wt\gh$ and $\pi:P\rightarrow M$ a $G$-structure  with $G\simeq G_\sharp/N_\sharp$. Then
\begin{itemize}
\item[i)] if $P$ is immersible to order $p$ in $\wt P$ but there exist $x_o\in M$ and $r\geq 0$ such that $[\widehat\Omega^{p-1}_r|_{x_o}]\neq 0$ for any admissible tuple $(s,\omega^0,\ldots,\omega^{p-1})$
then $P$ is not immersible to order $p+1$; 
\item[ii)] if $H^{p,2}(\wt\gh,W)=0$ for all $p\geq 0$ then $P$ is locally immersible in $\wt P$.
\end{itemize}
\end{corollary}
\section{Immersions into flat $2$-reductive $\wt G$-structures}\label{tortaele}
\setcounter{section}{4}
\setcounter{equation}{0}
\smallskip\par
We now restrict to the case of flat $2$-reductive $\wt G$-structures and $(\wt H/\wt K,\wt P)$ always denotes a $\wt G$-structure of this kind (hence with $\displaystyle\wt\gh=\wt\gh^{-1}+\wt\gh^0+\wt\gh^1$). We remark that any flat homogeneous $\wt G$-structure with $\wt G\subset \GL(V)$ irreducible on $V$ and of finite type is either $1$-reductive or $2$-reductive and that, in the second case, $\displaystyle\wt\gh$ is always determined by a grading of depth $1$ of a simple Lie algebra. For the classification of flat $2$-reductive $\wt G$-structures see \cite{KoNa, Oc1}.
\begin{definition}
\label{EQU}
Let $\pi:P\rightarrow M$ be a $G$-structure with $G\simeq G_\sharp/N_\sharp$.
Two admissible tuples $(s,\omega^0)$ and $(s',\omega'^0)$ on the same $\cU\subset M$ are {\it strongly equivalent} if $s=s'$ and $\omega'^0=\omega^0+[\omega^{-1},\varpi^{1}]$
for some $\wt\gh^{1}$-valued function $\varpi^{1}$.
\end{definition}
\par
This definition is motivated by the fact that the essential curvatures of two strongly equivalent tuples are necessarily the same. This can be checked as follows.
Let $(s,\omega^0)$, $(s,\omega'^0)$ be strongly equivalent.
By definitions, \eqref{MCC} and admissibility, one has $\omega^{-1}=\omega'^{-1}$, $\Omega^{-1}=\Omega'^{-1}$ and 
\beq
\label{febbre}
\Omega^{-1}=-[\omega^{-1},\omega^{0}]=-[\omega^{-1},\omega'^0]\ .
\eeq 
From this and strong equivalence it follows that
\begin{equation}
\begin{split}
\label{pranzo}
\Omega'^{0}&=d\omega'^0+\frac{1}{2}[\omega'^0,\omega'^0]\\
&=\Omega^0+d[\omega^{-1},\varpi^1]+[\omega^0,[\omega^{-1},\varpi^1]]+\frac{1}{2}[[\omega^{-1},\varpi^1],[\omega^{-1},\varpi^{1}]]\\
&=\Omega^0+[\Omega^{-1},\varpi^1]-[\omega^{-1},d\varpi^1]\\
&\qquad\;\,+[\omega^0,[\omega^{-1},\varpi^1]]+\frac{1}{2}[\omega^{-1},[\varpi^1,[\omega^{-1},\varpi^{1}]]\\
&\!\!\!\overset{\eqref{febbre}}=\Omega^0-[\omega^{-1},\epsilon^1]\\
&\qquad\;\;\text{with}\;\;\epsilon^1=d\varpi^1+[\omega^0,\varpi^1]+\frac{1}{2}[[\omega^{-1},\varpi^1],\varpi^1]\,,
\end{split}
\end{equation}
hence $[\Omega^{0}]=[\Omega'^{0}]$.
\smallskip\par
We also remark that if $H^{1,1}(\wt\gh,W)=0$ two tuples $(s,\omega^0)$ and $(s',\omega'^0)$ on the same contractible $\cU$ are strongly equivalent if and only if $s=s'$. Indeed if $s=s'$ equation \eqref{febbre} holds and $\omega'^0-\omega^0$ is closed, hence exact. 

A similar argument shows that when $H^{2,1}(\wt\gh,W)=0$ (that is $Z^{2,1}(\wt\gh,W)=0$ since $\wt\gh^2=0$) two admissible tuples $(s,\omega^0,\omega^1)$ and $(s',\omega'^0,\omega'^1)$ are equal if and only if $s=s'$ and $\omega^0=\omega'^0$.
\smallskip\par
In the following $\pi:P\rightarrow M$ is a $G$-structure with $G\simeq G_\sharp/N_\sharp$, $G_\sharp, N_\sharp\subset \wt G$ as in \eqref{semi1}, \eqref{semi2}. Furthermore, for any admissible tuple $(s,\omega^0)$ of $P$, we denote by 
$\cA_{(s,\omega^0)}$ the set of all admissible tuples $(s',\omega'^0,\omega'^1)$ with $s'=s$ and $\omega'^0=\omega^0$ and by $\cR$ the map which associates to any tuple in $\cA_{(s,\omega^0)}$ the corresponding total curvature,
\beq
\label{spiro}
\cR:\cA_{(s,\omega^0)}\to \cC^{\infty}(\cU, Z^{2,2}(\wt\gh,W))\;,\qquad\cR((s,\omega^0,\omega^1)):=\Omega^1\ .
\eeq
The following proposition shows that one can interchange any two strongly equivalent tuples in the analysis of the obstructions to the generalized integrability problem.
\begin{proposition}\label{finesez}
Let $(s,\omega^0)$ and $(s',\omega'^0)$ be two strongly equivalent tuples and $\cR$ and $\cR'$ the associated maps defined in \eqref{spiro}. Then 
$$
\cR(\cA_{(s,\omega^0)})=\cR'(\cA_{(s',\omega'^0)})\ .
$$ 
\end{proposition}
\begin{proof}
Let $\varpi^1$ be a $1$-form satisfying $\omega'^0=\omega^0+[\omega^{-1},\varpi^1]$ and $\epsilon^1$ as in \eqref{pranzo}. We observe that $\omega^1$ is a solution of $\Omega^{0}=-[\omega^{-1},\omega^1]$ if and only if $\omega'^1=\omega^1+\epsilon^1$ is a solution of $\Omega'^{0}=-[\omega'^{-1},\omega'^1]$.
This shows that the map $$
\cT:\cA_{(s,\omega^0)}\rightarrow\cA_{(s',\omega'^0)}\qquad \cT ((s,\omega^0,\omega^1)):=(s',\omega'^0,\omega^1+\epsilon^1)\ 
$$
is a bijection. Hence the claim is proved if we can show that the total curvature $\Omega'^1$ of $\cT(s,\omega^0,\omega^1)$ is equal to the total curvature $\Omega^1$ of $(s,\omega^0,\omega^1)$. By construction, standard properties of Lie algebra valued differential forms and the identities $\Omega^{-1}=-[\omega^{-1},\omega^0]$ and $\Omega^{0}=-[\omega^{-1},\omega^1]$
\begin{align*}
\Omega'^{1}&=d\omega'^1+[\omega'^0,\omega'^1]\\
&=\Omega^1+d\epsilon^1+[\omega^0,\epsilon^1]+[[\omega^{-1},\varpi^1],\omega^1]+[[\omega^{-1},\varpi^1],\epsilon^1]\\
&=\Omega^1+d[\omega^0,\varpi^1]+\frac{1}{2}d[[\omega^{-1},\varpi^1],\varpi^1]+[\omega^0,d\varpi^1]+[\omega^0,[\omega^0,\varpi^1]]\\
&\;\;\;\,+\frac{1}{2}[\omega^0,[[\omega^{-1},\varpi^1],\varpi^1]]+[[\omega^{-1},\varpi^1],\omega^1]
+[[\omega^{-1},\varpi^1],d\varpi^1]\\
&\;\;\;\,+[[\omega^{-1},\varpi^1],[\omega^0,\varpi^1]]+\frac{1}{2}[[\omega^{-1},\varpi^1],[[\omega^{-1},\varpi^1],\varpi^1]]\\
\end{align*}
\begin{align*}
\phantom{\Omega'^{1}}&=\Omega^1+\frac{1}{2}[[\Omega^{-1},\varpi^{1}],\varpi^1]
+\frac{1}{2}[[[\omega^{-1},\omega^{0}],\varpi^1],\varpi^1]+[\Omega^0,\varpi^1]\\
&\;\;\;\,+[[\omega^{-1},\omega^1],\varpi^1]-\frac{1}{2}[\omega^{-1},[[\omega^0,\varpi^1],\varpi^1]]-[\omega^{-1},[\omega^1,\varpi^1]]
\\
&\;\;\;\,-\frac{1}{2}[\omega^{-1},[d\varpi^1,\varpi^{1}]]-\frac{1}{6}[\omega^{-1},[[[\omega^{-1},\varpi^1],\varpi^1],\varpi^1]]\\
&=\Omega^1-\frac{1}{2}[\omega^{-1},[[\omega^0,\varpi^1],\varpi^1]]-[\omega^{-1},[\omega^1,\varpi^1]]
\\
&\;\;\;-\frac{1}{2}[\omega^{-1},[d\varpi^1,\varpi^{1}]]-\frac{1}{6}[\omega^{-1},[[[\omega^{-1},\varpi^1],\varpi^1],\varpi^1]]\ .
\end{align*}
Since $\wt\gh^2=0$, all terms except $\Omega^{1}$ are trivial, hence $\Omega^1=\Omega'^1$.  
\end{proof}
\begin{theorem}
\label{quaternion}
Let $H^{1,1}(\wt\gh,W)=H^{2,1}(\wt\gh,W)=0$. Given $0\leq p\leq 2$, if $P$ is locally immersible into $\wt P$ to order $p$, then:
\begin{itemize}
\item[i)] for any $x\in M$, the value $[\Omega^{p-1}|_x]$ of the essential $(G,\wt G)$-curvature of an admissible tuple $(s,\omega^0,\ldots,\omega^{p-1})$ depends only on the frame $u=s_x\in P|_x$; this determines a well-defined map
\beq
\label{curvatura}
\cR^{p+1}:P\rightarrow H^{p,2}(\wt\gh,W)\qquad \cR^{p+1}(u):=[\Omega^{p-1}|_x]\ ;
\eeq
\item[ii)]  the map \eqref{curvatura} is $G_\sharp$-equivariant with respect to its action  on $P$ via the isomorphism $G\simeq G_\sharp/N_\sharp$ and its standard actions on $W$ and $\wt\gh$;
\item[iii)] the map \eqref{curvatura} vanishes identically if and only if $P$ is locally immersible into $\wt P$ to order $p+1$.
\end{itemize}
\end{theorem}
\begin{proof}
Let $s$, $s':\cU\subset M\to P$ be two fixed local sections on the same contractible $\cU$, hence of the form
$s'=s\circ g^{-1}$ for some $g:\cU\rightarrow G$.
If $p_\sharp:G_\sharp\to G=G_\sharp/N_\sharp$ is the natural projection, then $g=p_\sharp\circ g_\sharp$ for some
$g_\sharp:\cU\rightarrow G_\sharp$.

As usual, for any admissible tuple $(s,\omega^0,\ldots,\omega^{p-1})$ associated with $s$, we use the frame $s_x$ to identify at each $x\in\cU$ the
$1$-forms $\omega^r|_x$ and the associated curvatures $\Omega^r$ with elements in $C(\wt\gh, W)$. On the other hand, for the admissible tuples $(s',\omega'^0,\ldots,\omega'^{p-1})$ associated with $s'$, we use $s'_x$ (and not $s_x$) to make the corresponding identifications for $\omega'^r|_x$ and $\Omega'^r$.
\vskip0.3cm\par\noindent
{\it Case $p=0$.} Let $\omega^{-1}$, $\omega'^{-1}$ be the $W$-valued $1$-forms associated with $s$, $s'$ as in \eqref{MCC} and $\Omega^{-1}$, $\Omega'^{-1}$ the corresponding total $(G,\wt G)$-curvatures. 
A direct computation shows that
\beq
\label{franca}
\Omega'^{-1}=g_\sharp\cdot\Omega^{-1}+\partial\left(R_{g_\sharp^{-1}}(dg_\sharp\circ s')\right)\ ,
\eeq
where we denote by $R_{g_\sharp}$ the natural right action of $g_\sharp$ on $TG_{\sharp}$. 
It follows that $[\Omega'^{-1}]=g_\sharp\cdot [\Omega^{-1}]$. 
From this, the transitivity of the action of $G$ on the fibers of $P$, and the fact that $s'_x=s_x$ implies $g_\sharp|_x=e$, (i) and (ii) follow. Claim (iii) follows from Theorem \ref{r-1}. 
\vskip0.3cm\par\noindent
{\it Case $p=1$.} 
Let $P$ be immersible in $\wt P$ to order $1$. Consider
an admissible tuple $(s,\omega^0)$ with associated total $(G,\wt G)$-curvature $\Omega^0$. By definition it satisfies $\Omega^{-1}=-\partial\omega^0$. Since $H^{1,1}(\wt\gh,W)=0$, we know that two admissible tuples $(s,\omega^0)$ and $(\widehat s,\widehat\omega^0)$ with $s=\widehat s$ are always strongly equivalent, hence with the same essential curvature of order $2$ (see remarks before Proposition \ref{finesez}). 
\vskip0.15cm\par
Now equation \eqref{franca} and 
the fact that the generalized Spencer operator is $G_\sharp$-equivariant yield
$$
\Omega'^{-1}=g_\sharp\cdot(-\partial\omega^0)+\partial\left(R_{g_\sharp^{-1}}(dg_\sharp\circ s')\right)=
-\partial\left(g_\sharp\cdot\omega^0-R_{g_\sharp^{-1}}(dg_\sharp\circ s')\right)\ ,
$$
that is $\left(s',g_\sharp{\cdot}(\omega^0{\circ} s^{-1})-R_{g_\sharp^{-1}}(dg_\sharp)\right)$ 
is admissible. Further, a direct computation shows that the total curvature $\Omega'^0$ of this admissible tuple is $\Omega'^0=g_\sharp\cdot\Omega^0$. By the above remark this implies that $[\Omega'^0]=g_\sharp\cdot[\Omega^0]$ independently of the choice of $\omega^0$ and $\omega'^0$. This yields (i), (ii) and (iii) with the same arguments of the case $p=0$. 
\vskip0.3cm\par\noindent
{\it Case $p=2$.} Let $P$ be immersible in $\wt P$ to order $2$. Since $H^{1,1}(\wt\gh,W)=H^{2,1}(\wt\gh,W)=0$, by the remarks before Proposition
\ref{finesez} any two admissible tuples $(s,\omega^0)$, $(\wh s,\wh \omega^0)$ with $s=\wh s$ are strongly equivalent and the sets $\cA_{(s,\omega^0)}$, $\cA_{(\wh s,\wh \omega^0)}$ consist both of just one element, say $\omega^1$ and $\wh\omega^1$ respectively. 
By Proposition \ref{finesez} 
the total $(G,\wt G)$-curvatures of $(s,\omega^0,\omega^1)$ and $(\wh s,\wh \omega^0,\wh \omega^1)$ coincide, meaning that such curvature depends only on $s$.
\vskip0.15cm\par
Let us consider now an admissible tuple $(s,\omega^0,\omega^1)$ and the admissible tuple $(s',\omega'^0,\omega'^1)$ where $\omega'^0$ is as in the previous case and $\omega'^1=g_\sharp{\cdot}(\omega^1{\circ}s^{-1})$.
The associated curvatures satisfy $\Omega'^1=g_\sharp\cdot\Omega^1$. 
By the above remark this implies that $\Omega'^1=g_\sharp\cdot\Omega^1$ independently of the choice of $(\omega^0,\omega^1)$ and $(\omega'^0,\omega'^1)$. This yields (i), (ii) and (iii) with the same arguments as above. 
\end{proof}
\vskip0.35cm\par
\section{Applications}
\label{teoria}
\setcounter{section}{5}
\setcounter{equation}{0}
\subsection{\it Riemannian immersions into spaces of constant curvature}
\label{Rie}
Let $(\wt M,\wt g)$ be a Riemannian manifold of dimension $\wt n$ and $\wt\pi:\wt P=O_{\wt g}(\wt M)\rightarrow \wt M$ the  $\OO_{\wt n}(\bR)$-structure given by the orthonormal frames 
$$
O_{\wt g}(\wt M)=\{(e_i)\in L(\wt M)\;|\; \wt g(e_i,e_j)=\delta_{ij}\;,\;\; 1\leq i,j\leq \wt n\}\ .
$$
Any $n$-dimensional submanifold $M\subset \wt M$ is $\wt P$-regular so that one can always consider the induced $G$-structure $\pi:P\rightarrow M$, which is the orthonormal frame bundle of the Riemannian metric $g$ induced on $M$ by $\wt g$.   
It follows that, given two Riemannian manifolds $(M,g)$ and $(\wt M,\wt g)$, the local immersibility of $\pi:O_g(M)\rightarrow M$ into $\wt\pi:O_{\wt g}(\wt M)\rightarrow\wt M$ is equivalent to the existence of local isometric immersions of $(M,g)$ into $(\wt M,\wt g)$.
\vskip0.3cm\par
Let us consider now the problem of the existence of a local isometric immersion of an $n$-dimensional Riemannian manifold $(M,g)$ into a space of constant curvature $k_o$, 
that is of the local immersibility of $P=O_g(M)$ in the homogeneous $\wt G$-structure $(\wt H/\wt K,\wt P)$ with
$$
\wt H=\left\{
\begin{array}{ll}
\SO_{\wt n}(\bR)\ltimes \bR^{\wt n}&\text{if}\;\;k_o=0\\
\SO_{\wt n+1}(\bR)&\text{if}\;\;k_o=1\\         
\SO_{\wt n,1}(\bR)&\text{if}\;\;k_o=-1
\end{array}\right.\;\;,\qquad\wt K=\SO_{\wt n}(\bR)\;\;,
$$
and $\wt P$ the orthonormal frame bundle of the homogeneous Riemannian manifold $(\wt M=\wt H/\wt K, \wt g)$, that is $\bR^{\wt n}$, $S^{\wt n}$ or $H^{\wt n}$ with their standard metrics of curvatures $k_o=0,1,-1$ respectively. 

These homogeneous spaces are reductive with quasi-graded Lie algebras
$$
\wt\gh=\wt\gh^{-1}+\wt\gh^0\;\;\;\;\;\text{where}\;\;\;\;\;\wt\gh^{-1}=\bR^{\wt n}\;\;,\;\;\;\wt\gh^0=\so_{\wt n}(\bR)\;\;,
$$
with Lie brackets defined by
$$
[A,v]=Av\in\wt\gh^{-1}\;\;,\;\;\;\;[v_1,v_2]=k_o(v_2\otimes \langle v_1,\cdot\rangle-v_1\otimes \langle v_2,\cdot\rangle)\in\wt\gh^0\ ,
$$
for any $v,v_1,v_2\in\wt\gh^{-1}$ and $A\in\wt\gh^0$. 
\vskip0.3cm\par
For our purposes, it is convenient to consider the decomposition of $\wt\gh^0$ into the following vector subspaces
\begin{align}
\nonumber
\so_{\wt n}(\bR)&{=}\left\{ \begin{pmatrix} A & -B^t \\ B & D \end{pmatrix}\,|\,A\in\so_n(\bR)\ ,\ B\in \bR^{\wt n-n}\otimes (\bR^{n})^*\ ,\ D\in\so_{\wt n-n}(\bR) \right\}\\
\label{serveafter}
&{=}\so(W)+W^\perp\otimes W^*+\so(W^\perp)\;\text{with}\;W=\bR^{n}\;,\; W^\perp=\bR^{\wt n-n}\,.
\end{align}
In $\wt\gh^0$ there is only one annihilator, that is
$\gc^0_1=\so(W^\perp)$, and from now on we consider $\gc^{0\perp}_1=\so(W)+W^\perp\otimes W^*$ as a fixed complementary subspace.
\vskip0.3cm\par
By the results of \S\ref{simdor} there exists a local isometric immersion of $(M,g)$ in $(\wt M,\wt g)$
if and only if there is a pair $(s,\omega^0)$ consisting of:
\begin{itemize}
\item[1)] a local section $s:\cU\subset M\rightarrow O_g(M)$, i.e. a field of $g$-orthonormal frames which we use to identify any tangent space $T_xM$ with $W=\bR^n$ and to define the $1$-form
$$\omega^{-1}|_x(w)=(w^1,\ldots,w^n,0,\ldots,0)\,,\;\;\;w=(w^i)\in W\simeq T_xM\,,\;\;\;x\in\cU\,;$$
\item[2)] a $\so_{\wt n}(\bR)$-valued local $1$-form $\omega^0$ satisfying the equations
\begin{align}
\label{muz}
&[\omega^{-1},\omega^0]=-d\omega^{-1}\ ,\\
\label{muzsecond}
&\Omega^0=d\omega^0+\frac{1}{2}[\omega^0,\omega^0]=0\ .
\end{align}
\end{itemize}
According to \eqref{serveafter} we decompose
$
\omega^0=\omega^{\so(W)}+\omega^{W^\perp\otimes W^*}+\omega^{\so(W^\perp)}
$ 
so that \eqref{muz} is equivalent to the following equations for any $x\in\cU$ and $w_1,w_2\in W$
\begin{align}
\label{levico}
&-\omega^{\so(W)}_x(w_2)(w_1)+\omega^{\so(W)}_x(w_1)(w_2)=-d\omega^{-1}|_{x}(w_1,w_2)\ ,\\
\label{levicoII}
&-\omega^{W^\perp\otimes W^*}_x(w_2)(w_1)+\omega^{W^\perp\otimes W^*}_x(w_1)(w_2)=0\ .
\end{align}
By classical arguments equation \eqref{levico} has a unique solution $\omega^{\so(W)}$, namely the connection $1$-form of the Levi-Civita covariant derivative $\nabla^{LC}$ of $(M,g)$.
On the other hand, if $NM$ is the trivial vector bundle $NM=\cU\times W^\perp\rightarrow \cU$, equation \eqref{levicoII} means that $\omega^{W^\perp\otimes W^*}$ can be identified with a symmetric local section $\Pi$ of $\otimes^2TM^*\otimes NM$.
We remark that, by the identification $W^\perp=\bR^{\wt n-n}$, the trivial bundle $NM$ has a natural fiber metric, which we denote by $g^\perp$, and $\omega^{\so(W^\perp)}$ is the connection $1$-form of a covariant derivative $\nabla^{\perp}$ on $NM$ compatible with $g^\perp$.
\smallskip\par
Using again \eqref{serveafter} we now decompose $\Omega^0$ into
$$
\Omega^0=\Omega^{\so(W)}+\Omega^{W^\perp\otimes W^*}+\Omega^{\so(W^\perp)}\ 
$$
and we observe that the curvature $\wh\Omega^0_1$ of level $1$ defined in \eqref{tavolo} and the complementary curvature $\wt\Omega^0_0$ of level $0$ defined in \eqref{tavolino} are equal to 
$$\wh\Omega^0_1=\Omega^{\so(W)}+\Omega^{W^\perp\otimes W^*}\;\;\;\;\text{and}\;\;\;\;\wt\Omega^0_0=\Omega^{\so(W^\perp)}\ .$$ 
Recall that $\Omega^0=\wh\Omega^0_1+\wt\Omega^0_0$. 
A simple computation shows
\begin{align*}
\Omega^{\so(W)}&=d\omega^{\so(W)}+\frac{1}{2}[\omega^{\so(W)},\omega^{\so(W)}]-(\omega^{W^\perp\otimes W^*})^{t}\cdot\omega^{W^\perp\otimes W^*}\\
&\phantom{=iiiiiiiiii,}+\frac{1}{2}[\omega^{-1},\omega^{-1}]\,,\\
\Omega^{W^\perp\otimes W^*}&=d\omega^{W^\perp\otimes W^*}+\omega^{W^\perp\otimes W^*}\cdot\omega^{\so(W)}+\omega^{\so(W^\perp)}\cdot\omega^{W^\perp\otimes W^*}\,,\\
\Omega^{\so(W^\perp)}&=d\omega^{\so(W^\perp)}+\frac{1}{2}[\omega^{\so(W^\perp)},\omega^{\so(W^\perp)}]
-\omega^{W^\perp\otimes W^*}\cdot(\omega^{W^\perp\otimes W^*})^t\ .
\end{align*}
Let $\nabla=\nabla^{LC}+\nabla^{\perp}$ be the covariant derivative on $TM\oplus NM$ determined by the connection $1$-form $\omega^{\so(W)}+\omega^{\so(W^\perp)}$ and denote by $R$ its curvature. We remark that $\nabla$ is compatible with the fiber metric $g+g^\perp$.

Equation \eqref{muzsecond} splits into the equations $\wh\Omega^0_1=\Omega^{\so(W)}+\Omega^{W^\perp\otimes W^*}=0$ and $\wh\Omega^0_0=\Omega^{\so(W^\perp)}=0$. One can check that the first corresponds to the Gauss-Codazzi equations
\vskip0.3cm\par\noindent
\begin{multline*}
g(R_{XY}Z,W)-g^\perp(\Pi(Y,Z),\Pi(X,W))+g^\perp(\Pi(X,Z),\Pi(Y,W))\\
+k_o(g(X,Z)g(Y,W)-g(Y,Z)g(X,W))=0\,,
\end{multline*}
and
$$
(\nabla_{X}\Pi)(Y,Z)-(\nabla_{Y}\Pi)(X,Z)=0\,,
$$
for any $X,Y,Z,W\in\mathfrak{X}(M)$. The second one corresponds to the Ricci equations
\begin{multline*}
g^\perp(R_{XY}\mu,\nu)+\sum_{i=1}^n g^\perp(\Pi(X,e_i),\mu)g^\perp(\Pi(Y,e_i),\nu)\\
-\sum_{i=1}^n g^\perp(\Pi(X,e_i),\nu)g^\perp(\Pi(Y,e_i),\mu)=0\ ,
\end{multline*}
for any $X,Y\in\mathfrak{X}(M)$ and local sections $\mu,\nu$ of $NM$, 
with $\{e_i\}$ a fixed local orthonormal frame field of $(M,g)$.
\par
By these observations it follows that, in this case, the results of \S\ref{simdor} yield the classical result that {\it a Riemannian manifold $(M,g)$ has a local isometric immersion in a space of constant curvature $k_o$ if and only if there exist a metric bundle $(NM,g^\perp)$, a compatible connection $\nabla^{\perp}$, and a symmetric tensor $\Pi$ that satisfy the Gauss-Codazzi and Ricci equations} (see e.g. \cite{Sp}). 
\medskip\par
\begin{remark}
\label{cisiamo!}
For Riemannian immersions into Euclidean spaces, i.e. with $k_o=0$, 
the homogeneous $\wt G$-structure $(\wt H/\wt K,\wt P)$ is flat $k$-reductive and one can discuss the obstructions to the generalized integrability problem as in \S\ref{biscotto}. 
In this situation the only relevant cohomology groups are $H^{0,2}(\wt\gg,W)$ and $H^{1,2}(\wt\gg,W)$ with $\wt\gg=\so_{\wt n}(\bR)$ and $W=\bR^{n}$. One can directly check that
$H^{0,2}(\so_{\wt n}(\bR),\bR^n)=0$ (this is equivalent to the fact that \eqref{levico} and \eqref{levicoII} are always solvable) and
$$
H^{1,2}(\so_{\wt n}(\bR),W)\simeq H^{1,2}(\so_{n}(\bR),W)+ W^ \perp\otimes \text{R}^{2,1}+\so(W^\perp)\otimes\Lambda^2 W^{*}\;,
$$ where R$^{2,1}$ is the unique $\ggl(W)$-irreducible submodule of $W^{*}\otimes\Lambda^2 W^{*}$ which is different from $\Lambda^3 W^{*}$ (that is the kernel of complete antysimmetrization).

As $\wt\gh^{1}=0$, the cohomology class of $\Omega^0$ is trivial
if and only if each of the three components
$$\Omega^{\mathfrak{so}(W)}|_{x}\in H^{1,2}(\so_{n}(\bR),W)\;\;,\;\;\;\;\Omega^{W^\perp\otimes W^*}|_{x}\in W^ \perp\otimes \text{R}^{2,1}\;\;\;\;\text{and}
$$
$$
\Omega^{\mathfrak{so}(W)}|_{x}\in \so(W^\perp)\otimes\Lambda^2 W^{*}\;\;,\phantom{cccccccccccccccccccccccccccccccccc}$$
is equal to zero, that is the Gauss, the Codazzi and the Ricci equations  are respectively satisfied.
\end{remark}
\begin{remark}\label{cisiamo!!}
We observe that the cohomology groups introduced in \S\ref{biscotto} are different from those considered by Rosly and Schwarz in \cite{RS}. In that paper the authors claimed that the only obstructions to the generalized integrability problem into flat $k$-reductive $\wt G$-structures are represented by the non-trivial classes in such cohomology groups (\cite[Appendix C]{RS}). On the other hand, in the case of Riemannian immersions into an Euclidean space one can check that the only non-trivial Rosly and Schwarz cohomology group is identifiable with the subspace $H^{1,2}(\so_n(\bR),\bR^n)$ of $H^{1,2}(\so_{\wt n}(\bR),\bR^n)$. This means that the set of Rosly and Schwarz obstructions correspond  just to the Gauss equation and seems therefore not to be complete.
\end{remark}
\subsection{\it Conformal immersions into a conformally flat space}
\label{Conf}\label{conforme}
Let $(\wt M,\wt g)$ be an $\wt n$-dimensional Riemannian manifold and $\wt\pi:\wt P=\CO_{\wt g}(\wt M)\rightarrow \wt M$ the $\CO_{\wt n}(\bR)$-structure of conformal frames 
$$
\CO_{\wt g}(\wt M)=\left\{(e_i)\in L(\wt M)\;|\;\wt g(e_i,e_j)=c\delta_{ij}\;,\;1\leq i,j\leq \wt n\;,\;c\in(0,+\infty)\right\}\,.
$$
Any $n$-dimensional submanifold $M\subset \wt M$ is $\wt P$-regular and has an induced $G$-structure $\pi:P\rightarrow M$, naturally identifiable with the bundle of conformal frames of the induced Riemannian metric $g$. It follows that, given two Riemannian manifolds $(M,g)$, $(\wt M,\wt g)$, the existence of local conformal immersions of $(M,g)$ into $(\wt M,\wt g)$ is equivalent to the existence of a local immersion of $\pi:\CO_g(M)\rightarrow M$ into $\wt\pi:\CO_{\wt g}(\wt M)\rightarrow\wt M$.
Due to this, the existence of a local conformal immersion of $(M,g)$ in the homogeneous conformally flat space $(S^{\wt n},\wt g)$ is equivalent to the local immersibility of $\pi:\CO_g(M)\rightarrow M$ in the homogeneous $\wt G$-structure $(S^{\wt n}=\SO_{\wt n+1,1}(\bR)/\wt K,\wt P)$ with $\wt P=\CO_{\wt g}(S^{\wt n})$ and $\wt K$ an appropriate parabolic subgroup.

This homogeneous $\wt G$-structure is flat $k$-reductive and associated with the graded Lie algebra $\wt\gh=\mathfrak{so}_{\wt{n}+1,1}(\bR)$ with grading
$$
\wt\gh=\wt\gh^{-1}+\wt\gh^{0}+\wt\gh^{1}= V+\mathfrak{co}(V)+V^*\;\;,\;\;\;\text{where}\;\;\;\;\;V=\bR^{\wt n}\;\;,
$$
and Lie brackets determined by the standard action of $\mathfrak{co}(V)$ on $V$ and $V^*$ and the bracket between $\alpha\in V^*$ and $v\in V$ defined by
\beq
\label{ancoraa}
[\a,v]=v\otimes \a-\natural(\a\otimes v)+\a(v)I\,.
\eeq
Here $\natural:V^*\otimes V\rightarrow V\otimes V^*$ is the natural isomorphism induced by the standard scalar product $\langle\cdot,\cdot\rangle$ of $V=\bR^{\wt n}$. The isotropy subalgebra
is $$\wt\gk=\wt\gh^0+\wt\gh^1=\mathfrak{co}(V)+V^*\ .$$
We recall that when $\wt n\geq 3$, $\wt\gh$ is the maximal prolongation of $\mathfrak{co}(V)$
(\cite{Kb, SS}). 
\smallskip\par
For the case of conformal immersions into a conformally flat Riemannian space,
the analogues of the Gauss, Codazzi and Ricci equations are provided by the following theorem.

We first need to fix some notation. Let $(M,g)$ be an $n$-dimensional Riemannian manifold with Levi-Civita connection $\nabla^{LC}$
and $\wt n\geq n$. Given two symmetric tensors $b$, $b'$ of type $(0,2)$ on $M$, we also denote by $b\varowedge b'$ their Kulkarni-Nomizu product, namely the $(0,4)$-tensor defined by
\begin{multline}\nonumber
b\varowedge b'(X,Y,Z,W)=b(X,Z)b'(Y,W)+b(Y,W)b'(X,Z)\\
-b(X,W)b'(Y,Z)-b(Y,Z)b'(X,W)\,.
\end{multline}
\begin{theorem}
\label{miglioverdeprima}
Given a conformal immersion $\imath:(M,g)\longrightarrow(\bR^{\wt n},\langle\cdot,\cdot\rangle)$
with $\imath^*\langle\cdot,\cdot\rangle=e^{2f}g$ for a smooth function $f$, the quadruple 
$(\nabla^{\perp},\Pi, B, D)$ formed by the restriction $\nabla^\perp$  of the Levi-Civita connection of $(\bR^{\wt n},\langle\cdot,\cdot\rangle)$ along the normal bundle $NM$ of $M$ and the tensors
\begin{align*}
\Pi&=e^{-f}\cdot\overline\Pi\,,\\
B&=\operatorname{Hess}(f)-df\odot df+\frac{1}{2}g(\nabla f,\nabla f)g\,,\\
D&=e^{f}\cdot\overline\Pi(\nabla f,\cdot)\,,
\end{align*}
where $\overline\Pi$ is the second fundamental form of $(M,e^{2f}g)\subset(\bR^{\wt n},\langle\cdot,\cdot\rangle)$ and $\nabla f$ and $\operatorname{Hess}(f)$ are the gradient and the Hessian of $f$ w.r.t. $g$, satisfies the following equations for any $X,Y,Z,W{\in}\mathfrak{X}(M)$ and $\mu,\nu{\in}\Gamma(NM)$:
\begin{align}
\nonumber
&g(R_{XY}Z,W)-\langle\Pi(Y,Z),\Pi(X,W)\rangle+\langle\Pi(X,Z),\Pi(Y,W)\rangle\\
\label{la1}
&\qquad\qquad\qquad=-g\varowedge B(X,Y,Z,W)\,,\\
\label{la2}
&(\nabla_{X}\Pi)(Y,Z)-(\nabla_{Y}\Pi)(X,Z)=D(X)g(Y,Z)-D(Y)g(X,Z)\,,\\
\nonumber
&\langle R_{XY}\mu,\nu\rangle+\sum_{i=1}^n \langle\Pi(X,e_i),\mu\rangle \langle\Pi(Y,e_i),\nu\rangle\\
\label{la3}
&\qquad\qquad\qquad-\sum_{i=1}^n \langle\Pi(X,e_i),\nu\rangle \langle\Pi(Y,e_i),\mu\rangle=0\,,\\
\label{la4}
&(\nabla_{X}B)(Y,Z)+\langle D(X),\Pi(Y,Z)\rangle\;\text{is symmetric in}\; X, Y,\\
\label{la5}
&(\nabla_{X}D_\flat)(Y,\mu)+\sum_{i=1}^n \langle\Pi(X,E_i),\mu\rangle B(Y,E_i)\;\text{is symmetric in}\; X, Y,
\end{align}
where $R$ is the curvature of the covariant derivative $\nabla=\nabla^{LC}+\nabla^{\perp}$ on the vector bundle $TM\oplus NM\to M$, $\{e_i\}$ a fixed local orthonormal frame field of $(M,g)$ and
$\flat:NM\rightarrow NM^*$ the natural isomorphism induced by $\langle \cdot,\cdot\rangle|_{NM}$.
\end{theorem}
\begin{proof}
The equations \eqref{la1}-\eqref{la3} are a consequence of the Gauss-Codazzi-Ricci equations for the isometric immersion $(M,e^{2f}g)\subset (\bR^{\wt n},\langle\cdot,\cdot\rangle)$ and the explicit expressions of Levi-Civita connection and curvature of $g'=e^{2f}g$ in terms of $g$ and $f$ (see e.g. \cite[Thm. 1.159]{B}; caution: we use a definition of curvature which is opposite in sign to the one in \cite{B}). Equations \eqref{la4} and \eqref{la5} follow from the definitions of $B$ and $D$ by a somehow long but straightforward computation.
\end{proof}
The contents of \S\ref{biscotto} yield the following theorem which shows that the above equations are actually necessary and sufficient conditions for the existence of local conformal immersions into conformally flat spaces. This is a result that was actually first obtained by Akivis (\cite{AK61}, see also \cite[Theorem 3.1.5]{AG}) using the method of moving frames. 
\par\smallskip
In the statement of this theorem, given an open subset $\cU\subset M$, we denote by $NM=\cU\times W^\perp\rightarrow \cU$ the trivial vector bundle with fiber $W^\perp=\bR^{\wt n-n}$ and natural fiber metric
determined by the standard inner product $\langle \cdot,\cdot\rangle|_{W^\perp}$.
\begin{theorem}
\label{miglioverde}
{\it 
Let $(M,g)$ be an $n$-dimensional Riemannian manifold and $\wt n \geq n$, $\wt n\geq 3$. Then, for any $x\in M$, there exists a local conformal immersion of $(M,g)$ around $x$ into $(S^{\wt n},\wt g)$ if and only if for some neighborhood $\cU\subset M$ of $x$ there is a quadruple $(\nabla^\perp,\Pi,B,D)$ consisting of} 
\begin{itemize}
\item[i)] {\it a metric connection $\nabla^\perp:TM|_{\cU}\times NM\rightarrow NM$,}
\item[ii)] {\it a local section $\Pi$ of $S^2T^*M\otimes NM$,}
\item[iii)] {\it a local section $B$ of $S^2T^*M$,}
\item[iv)] {\it a local section $D$ of $T^*M\otimes NM$,}
\end{itemize}
which satisfies \eqref{la1}$-$\eqref{la5}.
\end{theorem}
\begin{proof}
According to the decomposition \eqref{serveafter} of $\so(V)$ we consider the decomposition of $\wt\gh^0=\mathfrak{co}(V)$ into the vector subspaces  
\beq
\label{serveaftercon}
\mathfrak{co}(V)=\so(V)+\bR=\so(W)+W^\perp\otimes W^*+\so(W^\perp)+\bR\ .
\eeq
We also consider the decomposition of $\wt\gh^{1}=V^*$ given by
\beq
\label{ancoraa2}
V^*=W^*+(W^{\perp})^*\ .
\eeq
Since $\wt\gh^1\neq (0)$, by the results of \S\ref{simdor} there exists a local conformal immersion of $(M,g)$ into $(S^{\wt n},\wt g)$ if and only if there is an admissible triple $(s,\omega^0,\omega^1)$, that is a triple consisting of:
\begin{itemize}
\item[1)] a local section $s:\cU\subset M\rightarrow \CO_g(M)$, i.e. a field of $g$-conformal frames which we use to identify any tangent space $T_xM$ with $W=\bR^n$ and define the $1$-form
$$\omega^{-1}|_x(w)=(w^1,\ldots,w^n,0,\ldots,0)\,,\;\;\;w=(w^i)\in W\simeq T_xM\,,\;\;\;x\in\cU\ .$$
By the observation before Definition \ref{acqua}  we may assume without loss of generality that $s$ is a field of $g$-orthonormal frames (in particular   $g$ is identified at any point with the standard scalar product $\langle\cdot,\cdot\rangle|_{W}$);
\item[2)] a $\mathfrak{co}(V)$-valued local $1$-form $\omega^0$ satisfying the equation
\begin{align}
\label{muz2}
&[\omega^{-1},\omega^0]=-d\omega^{-1};
\end{align}
\item[3)] a $V^*$-valued local $1$-form $\omega^1$ satisfying the equations
\begin{align}
\label{pizzastas}
&[\omega^{-1},\omega^1]=-\Omega^0=-d\omega^0-\frac{1}{2}[\omega^0,\omega^0]\;,\\
\label{pizzastas2}
&\Omega^1=d\omega^1+[\omega^1,\omega^0]=0\;.
\end{align}
\end{itemize}
\medskip\par
According to \eqref{serveaftercon} we decompose
$
\omega^0=\omega^{\so(W)}+\omega^{W^\perp\otimes W^*}+\omega^{\so(W^\perp)}+\omega^{\bR}
$ 
and \eqref{muz2} is equivalent to the following equations for any $x\in\cU$, $w_1,w_2\in W$:
\begin{align}
\nonumber
-d\omega^{-1}|_{x}(w_1,w_2)=&-\omega^{\so(W)}_x(w_2)(w_1)+\omega^{\so(W)}_x(w_1)(w_2)\\
\label{simba90}
&\qquad\qquad-\omega_x^{\bR}(w_2)(w_1)+\omega_x^\bR(w_1)(w_2)\;,\\
\label{simba91}
0=&-\omega^{W^\perp\otimes W^*}_x(w_2)(w_1)+\omega^{W^\perp\otimes W^*}_x(w_1)(w_2)\;.
\end{align}
From now on we consider the solution of \eqref{simba90} determined by $\omega^\bR=0$ and the connection $1$-form $\omega^{\so(W)}$ of $\nabla^{LC}$. Equation \eqref{simba91}  means that $\omega^{W^\perp\otimes W^*}$ can be identified with a symmetric local section $\Pi$ of $\otimes^2TM^*\otimes NM$. Finally we observe that $\omega^{\so(W^\perp)}$ is the connection $1$-form of a metric connection $\nabla^{\perp}$ on $NM$.
\medskip\par
Since $\omega^\bR=0$ we can argue in complete analogy with \S\ref{Rie} and decompose $\Omega^0$ according to \eqref{serveaftercon}
$$
\Omega^0=\Omega^{\so(W)}+\Omega^{W^\perp\otimes W^*}+\Omega^{\so(W^\perp)}\ .
$$
The explicit expressions of each component
$\Omega^{\so(W)}$, $\Omega^{W^\perp\otimes W^*}$ and $\Omega^{\so(W^\perp)}$ is as in \S\ref{Rie}.
Further, according to \eqref{ancoraa2}, we decompose
$
\omega^1=\omega^{W^*}+\omega^{(W^\perp)^*}
$
so that equation \eqref{pizzastas} splits into the following set of equations
\begin{align*}
\label{tttt}
\Omega^{\so(W)}(w_1,w_2)&=-w_2\otimes\omega^{W^*}(w_1)+\natural(\omega^{W^*}(w_1)\otimes w_2)\\
&\phantom{=\,}+w_1\otimes \omega^{W^*}(w_2)-\natural(\omega^{W^*}(w_2)\otimes w_1)\,,\\
\Omega^{W^\perp\otimes W^*}(w_1,w_2)&=\natural(\omega^{(W^\perp)^*}(w_1)\otimes w_2)-\natural(\omega^{(W^\perp)^*}(w_2)\otimes w_1)\,,\\
\Omega^{\so(W^\perp)}(w_1,w_2)&=0\,,\\
\omega^{W^*}(w_1)(w_2)&=\omega^{W^*}(w_2)(w_1)\ .
\end{align*}
If we denote by $B$ (resp. $D$) the local section of $\otimes^2 T^*M$ (resp. $T^*M\otimes NM$)
which corresponds to the local $1$-form $\omega^{W^*}$ (resp. $\omega^{(W^\perp)^*}$), one can check that the above equations correspond to \eqref{la1}, \eqref{la2}, \eqref{la3} and the fact that $B$ is symmetric.
\vskip0.3cm\par
Now according to \eqref{ancoraa2} we decompose $\Omega^1$ into
$\Omega^1=\Omega^{W^*}+\Omega^{(W^\perp)^*}$.
A simple computation shows 
\begin{align*}
\Omega^{W^*}&=d\omega^{W^*}-\omega^{W^*}\cdot \omega^{\so(W)}-\omega^{(W^\perp)^*}\cdot\omega^{W^\perp\otimes W^*}\,,\\
\Omega^{(W^\perp)^*}&=d\omega^{(W^\perp)^*}-\omega^{W^*}\cdot(\omega^{W^\perp\otimes W^*})^t
-\omega^{(W^\perp)^*}\cdot\omega^{\mathfrak{so}(W^\perp)}\ .
\end{align*}
Equation \eqref{pizzastas2} splits into $\Omega^{W^*}=0$ and $\Omega^{(W^\perp)^*}=0$.
Using equations \eqref{la1}-\eqref{la3} one can check that the first equation correspond to \eqref{la4} and the second one to \eqref{la5}.
By these observations and the general results of \S\ref{simdor} it follows that a Riemannian manifold $(M,g)$ with a triple $(\nabla^{\perp}, \Pi, B, D)$ as in i)-iv) and satisfying \eqref{la1}-\eqref{la5} has a local conformal immersion in $(S^{\wt n},\wt g)$. This concludes the proof of the ``if'' direction of the Theorem. 
\medskip\par
The ``only if'' direction follows from Theorem \ref{miglioverdeprima} and the fact that $S^{\wt n}$ is locally conformal to $\bR^{\wt n}$.
\end{proof}
%
The results in \S\ref{biscotto} and \S\ref{tortaele} allow to improve the above theorem and show that in many cases the system \eqref{la1}-\eqref{la5} is equivalent to a smaller one.
\begin{theorem}
\label{elenaarriva}
Under the assumption of Theorem \ref{miglioverde} we have that:
\medskip\par\noindent
i) If $n\geq 4$ any solution $(\nabla^\perp,\Pi,B, D)$ of \eqref{la1}, \eqref{la2} and \eqref{la3} automatically satisfies \eqref{la4} and \eqref{la5} as well;
\medskip\par\noindent
ii) If $\wt n=n+1$, local conformal immersions as hypersurfaces exist if and only if there are
tensors $\Pi, B, D$ as in Theorem \ref{miglioverde} satisfying \eqref{la1}, \eqref{la2}, \eqref{la4} and \eqref{la5} with $\nabla^\perp$ the flat connection. Moreover
\begin{itemize}
\item[a)] if $n=2$ there is always a solution $(\Pi,B, D)$ of \eqref{la1} and \eqref{la2},
\item[b)] if $n\geq 4$ any solution $(\Pi,B, D)$ of \eqref{la1}, \eqref{la2} automatically satisfy \eqref{la4} and \eqref{la5} as well.
\end{itemize}
\end{theorem}
\begin{proof}
We first determine the cohomology groups $H^{p,2}(\wt\gh,W)$ associated with $\wt\gh=\so_{\wt n+1,1}(\bR)$ and $W=\bR^{n}$, for any $n\leq \wt n$. As $\wt\gh^2=0$ the only relevant cases are $p=0,1$ and $2$.
We observe that if $n=\wt n$ these groups coincide with the usual
Spencer groups $H^{p,2}(\mathfrak{co}_{n}(\bR))$ of $\mathfrak{co}_{n}(\bR)$ and
(see \cite{Oc2})
$$
H^{p,2}(\mathfrak{co}_{n}(\bR))\neq 0\;\;\;\;\text{only when}\;\;\;\;\;i)\;p=2,\;n=3\;\;\;\;\;\text{and}\;\;\;\;\;ii)\;p=1,\;n\geq 4\,.
$$
Further, for any $n<\wt n$, one can directly check that
$H^{0,2}(\wt\gh,W)=0$ (this is equivalent to the fact that \eqref{simba90} and \eqref{simba91} are always solvable).
\medskip\par
Consider an element $\omega^1\in C^{2,2}(\wt\gh,W)=V^*\otimes\Lambda^2 W^*$ and decompose it according to \eqref{ancoraa2} into
$
\omega^1=\omega^{W^*}+\omega^{(W^\perp)^*}$.
From definitions it follows
$$\partial\omega^{W^*}\in\mathfrak{co}(W)\otimes\Lambda^3 W^*\;\;,\;\;\;\;\;\;\partial\omega^{(W^\perp)^*}\in W^\perp\otimes W^*\otimes \Lambda^3 W^*\;\;,$$
so that the vanishing of $\partial\omega^1=\partial\omega^{W^*}+\partial\omega^{(W^\perp)^*}$ is equivalent to
\beq
\label{paura11}
\partial\omega^{W^*}=0
\eeq
and
\beq
\label{paura21}
\!\!\!\partial\omega^{(W^\perp)^*}=0\ .
\eeq
If $n=2$, $\Lambda^3 W^*=0$ so that \eqref{paura11} and \eqref{paura21} are always trivially satisfied.
If $n\geq 3$ one can directly show that \eqref{paura21} implies $\omega^{(W^\perp)^*}=0$.

As $\wt\gh^2=0$ one has $H^{2,2}(\wt\gh,W)=Z^{2,2}(\wt\gh,W)$ and by the above discussion
\beq
\label{echi1}
H^{2,2}(\wt\gh,W)\simeq
\begin{cases}
H^{2,2}(\mathfrak{co}_{n}(\bR))\;\;\;\text{if}\;\;\;n\geq 3\,,\\
V^*\otimes\Lambda^2W^*\;\;\;\;\,\text{if}\;\;\;n=2\,.
\end{cases}
\eeq
In particular $H^{2,2}(\wt\gh,W)$ is trivial when $n\geq 4$, for any $\wt n$.
\medskip\par
Consider now an element $\omega^0\in C^{1,2}(\wt\gh,W)=\mathfrak{co}(V)\otimes\Lambda^2 W^*$ and decompose it according to \eqref{serveaftercon} into 
$
\omega^0=\omega^{\so(W)}+\omega^{W^\perp\otimes W^*}+\omega^{\so(W^\perp)}+\omega^{\bR}
$. 
Equation $\partial\omega^0=0$ splits into the following two equations
\begin{align*}
\partial\omega^{\so(W)}+\partial\omega^{\bR}=0\qquad\text{and}\qquad\partial\omega^{W^\perp\otimes W^*}=0\ .
\end{align*}
Let also $\omega^1\in C^{2,1}(\wt\gh,W)=V^*\otimes W^*$ and decompose it according to \eqref{ancoraa2} 
into $\omega^1=\omega^{W^*}+\omega^{(W^\perp)^*}$.
Using \eqref{ancoraa} one can check that the equation $\partial\omega^1=\omega^0$ is equivalent to
$$
\partial\omega^{W^*}=\omega^{\so(W)}+\omega^{\bR}\;\;\;\;\;\;\;\;\text{and}\;\;\;\;\;\;\;\;\partial\omega^{(W^\perp)^*}=\omega^{W^\perp\otimes W^*}\ .
$$
The above discussion together with a direct computation and a dimensional argument
shows that
\begin{align}
\label{echi20}
H^{1,2}(\wt\gh,W)\simeq 
\begin{cases} 
H^{1,2}(\mathfrak{co}_{n}(\bR))+W^\perp\otimes\frac{\text{R}^{2,1}}{W^*}+\mathfrak{so}(W^\perp)\otimes \Lambda^2 W^*\;\;\text{if}\;\;\;n\geq 3\,,\\
W^\perp\otimes\frac{\text{R}^{2,1}}{W^*}+\mathfrak{so}(W^\perp)\otimes\Lambda^2 W^*\;\;\;\;\;\;\;\;\;\;\;\;\;\;\;\;\;\;\;\;\;\;\;\;\;\;\,\text{if}\;\;\;n=2\,,
\end{cases}
\end{align}
and $$\dim\frac{\text{R}^{2,1}}{W^*}=\frac{n^3-4n}{3}\,,$$ 
where R$^{2,1}\subset W^{*}\otimes\Lambda^2 W^{*}$ is the kernel of complete antysimmetrization.
In particular $H^{1,2}(\wt\gh,W)$ is trivial when $n=2$, $\wt n=3$. 
\medskip\par
Point (i) of the theorem is a direct consequence of \eqref{echi1}, Theorem \ref{r-1}, the proof of Theorem \ref{miglioverde} and, finally, of the observation that the essential curvature $[\Omega^1]$ is naturally identified with the total curvature $\Omega^1$. 
Point (ii) of the theorem follows similarly from \eqref{echi1} and \eqref{echi20} and by applying Theorem \ref{r-1} and Theorem \ref{miglioverde} when $\wt n=n+1$, $W^\perp$ is $1$-dimensional and $\so(W^\perp)=0$. 
\end{proof}
\begin{remark}
If $n=3$ neither the group $H^{1,2}(\wt\gh,W)$ nor $H^{2,2}(\wt\gh,W)$ vanishes and the system \eqref{la1}-\eqref{la5} can not be reduced simply as in Theorem \ref{elenaarriva}.
\end{remark}
 \subsection{\it CR immersions into $\bC^n$}\label{CRC}
Let $(\wt M,\wt J)$ be an almost complex manifold, a real manifold of even dimension $\wt n=2\wt m$ with a tensor $\wt J:T\wt M\rightarrow T\wt M$ satisfying $\wt J^2=-I$. Fix an integer $0\leq k\leq \wt m-1$ and consider 
the $\wt G$-structure on $\wt M$, $\wt G\simeq \GL_{\wt m}(\mathbb{C})$,
$$
\wt P_{\wt J}=\left\{(e_i)\in L(\wt M)\;|\;e_{\wt m-k+i}=\wt J e_i\;\;\text{for any}\;\; 1\leq i\leq \wt m-k\;\;\text{and}\;\;\;\;\;\;\;\;\;\;\;\;\;\;\;\;\;\;\;\;\;\;\;\;\;\right.
$$
$$
\;\;\;\;\;\;\;\qquad\qquad\qquad\left.\phantom{\wt M} e_{2\wt m-k+i}=\wt J e_{2(\wt m-k)+i}\;\;\text{for any}\;\; 1\leq i\leq k\right\}\,.
$$
One may check that a generic submanifold $M\subset \wt M$ of dimension $n=\wt n-k$ is $\wt P$-regular and it is endowed with
\begin{itemize}
\item[i)] a rank $2(\wt m-k)$ distribution $\cD\subset TM$,
\item[ii)] a tensor $J:\cD\rightarrow\cD$ satisfying $J^2=-I$,
\end{itemize} 
that is $(M,\cD,J)$ is an almost CR manifold of CR codimension $k$. 
The induced $G$-structure on $M$ is 
$$
P_{\cD,J}=\left\{(e_i)\in L(M)\phantom{\wt i}|\,e_{i}\in\cD\;\;\text{if}\;\;1\leq i\leq 2(\wt m-k)\;\;\text{and}
\;\;\;\;\;\;\;\;\;\;\;\;\;\;\;\;\;\;\;\;\;\;\;\;\;\right.
$$
$$
\;\;\;\;\;\;\;\;\;\left.\phantom{M}\;e_{\wt m-k+j}=J e_{j}\;\;\text{if}\;\; 1\leq j\leq \wt m-k\right\}\,.
$$
In this case the group $N_\sharp$ defined in \eqref{semi2} is trivial so that 
$G$ is naturally identified with the group $G_\sharp$ defined in \eqref{semi1}. 
\medskip\par
It follows that, given an almost CR manifold $(M,\cD,J)$ of real dimension $n=\wt n-k$ and CR codimension $k$ and an almost complex manifold $(\wt M,\wt J)$ of real dimension $\wt n=2\wt m$, the existence of a local CR immersion of $(M,\cD,J)$ into $(\wt M,\wt J)$ is equivalent to the existence of a local immersion of the associated bundle $P_{\cD,J}$ into the bundle $\wt P_{\wt J}$ of $\wt M$.
In particular, the existence of a local CR immersion of $(M,\cD,J)$ into $\bC^{\wt m}$, endowed with the standard complex structure $\wt J_{o}$, is equivalent to the local immersibility of $\pi:P_{\cD,J}\rightarrow M$ into the $\wt G$-structure $\wt\pi:\wt P_{\wt J_o}\to \bC^{\wt m}$. 

The latter is homogeneous and flat but it is not of finite type, since the maximal prolongation $\textstyle\wt\gh_\infty=\sum_{p= -1}^{\infty}\wt\gh_\infty^p$ of $\wt\gg=\ggl_{\wt m}(\bC)$ is infinite dimensional (see e.g. \cite{Kb}). However the notions of immersibility to order $p$ and of generalized cohomology still make sense and, in the real analytic setting, immersibility at any order still yields complete immersibility (\cite{G, RS}).
\medskip\par
The next proposition provides the relevant cohomology groups. We first need to fix some notation.
Let $W\simeq \bR^{n}$ (resp. $U\simeq\bR^{2(\wt m-k)}$) be the subspace of $V=\bR^{\wt n}$ determined by the vanishing of the last $k$ (resp. $2k$) standard coordinates. We also denote by $W^\perp$ (resp. $U^\perp$) the subspace of $V$ (resp. $W$) determined by the vanishing of the first $n$ (resp. $2(\wt m-k)$) standard coordinates of $V$. There are vector space direct sum decompositions 
$$
V=W+W^\perp\simeq \bR^{n}+\bR^{k}\;\;,\qquad W=U+U^\perp\simeq \bR^{2(\wt m-k)}+\bR^{k}\;.
$$
\begin{proposition}
\label{stef}
All groups $H^{p,2}(\wt\gh_\infty,W)$, $p\geq 1$, are trivial. 
\end{proposition}
\begin{proof}
We recall that the $(p-1)^{\text{th}}$-component of the maximal prolongation $\wt\gh_\infty$ of $\wt\gg=\ggl_{\wt m}(\bC)$ is 
$
\wt\gh_\infty^{p-1}=\bC^{\wt m}\otimes_{\bC} S^{p}(\bC^{\wt m})^* 
$
and it has a natural structure of complex vector space (see e.g. \cite{Kb}). 

We also note that, since any $v\in W^\perp$ satisfies $\wt J_o v\in U^\perp\subset W$, it is possible to extend any element
$
\omega^{p-1}\in C^{p,2}(\wt\gh_\infty,W)=\wt\gh_\infty^{p-1}\otimes\Lambda^2 W^*
$
to an element
$\wt\omega^{p-1}\in C^{p,2}(\wt\gh_\infty,V)=\wt\gh_\infty^{p-1}\otimes\Lambda^2 V^*$, 
setting
\begin{align*}
\wt\omega^{p-1}(v_1,v_2)&:=-\omega^{p-1}(\wt J_o v_1,\wt J_o v_2)&\qquad&\text{for all}\;\;\;v_1,v_2\in W^\perp\,,\\
\wt\omega^{p-1}(v_1,w_2)&:=-i\omega^{p-1}(\wt J_o v_1,w_2)&\qquad&\text{for all}\;\;\;v_1\in W^\perp\,,\;w_2\in W\,.
\end{align*}
We finally observe that the generalized Spencer operator on $C^{p,2}(\wt\gh_\infty,V)$ is just the usual Spencer operator and that the corresponding cohomology $H^{p,2}(\wt\gh_\infty,V)$ is just the Spencer group $H^{p,2}(\ggl_{\wt m}(\bC))$ of $\ggl_{\wt m}(\bC)$. It is known that this group is trivial for any $p\geq 1$ (\cite{G, KS}).  

Now, one can check that if $\partial\omega^{p-1}=0$ then also $\partial\wt\omega^{p-1}=0$. From this fact, for all $p\geq 1$ and $\omega^{p-1}$ with $\partial\omega^{p-1}=0$, one has
\begin{align*}\wt\omega^{p-1}&=\partial\wt\omega^{p}&\;\;\;\;\;\;&\text{for some}&\;\;\;\;\;\;&\wt\omega^p\in C^{p+1,1}(\wt\gh_\infty,V)\;\;\;\;\;\;\;\,\text{so that}\\
\omega^{p-1}&=\partial\omega^{p}&\;\;\;\;\;\;&\text{where}&\;\;\;\;\;\;&\omega^{p}=\wt\omega^{p}|_W\in C^{p+1,1}(\wt\gh_\infty,W)\ .
\end{align*}
In other words $H^{p,2}(\wt\gh_\infty,W)=0$ for every $p\geq 1$.
\end{proof}
Now, we observe that any coboundary $T\in B^{0,2}(\wt\gh_\infty,W)$ is the restriction to $\Lambda^2W$ of some $\wt T\in B^{0,2}(\wt\gh_\infty,V)$ and
$$B^{0,2}(\wt\gh_\infty,V)=\left\{\wt T\in V\otimes\Lambda^{2} V^*\;\;|\;\;\text{for every}\;\;v_1,v_2\in V\right.\qquad\qquad\qquad\qquad\qquad\;\;\;
$$
$$
\qquad\qquad\qquad\;\;\;\left.\wt T(v_1,v_2)-\wt T(\wt J_o v_1,\wt J_o v_2)=-\wt J_o\wt T(\wt J_ov_1,v_2)-\wt J_o\wt T(v_1,\wt J_o v_2)\right\}.
$$
From this, one gets that the vector space
$$
B^{0,2}(\wt\gh_\infty,W)\cap W\otimes\Lambda^2 W^*
$$
is formed by all maps $T\in W\otimes\Lambda^2W^*$ satisfying
\beq
\label{thee}
\begin{split}
T(u_1,u_2)-T(\wt J_o u_1,\wt J_o u_2)&\in U\;,\\
T(u_1,u_2)-T(\wt J_o u_1,\wt J_o u_2)&=-\wt J_oT(\wt J_o u_1,u_2)-\wt J_o T(u_1,\wt J_o u_2)\;,
\end{split}
\eeq
for any $u_1,u_2\in U$.
By Proposition \ref{stef}, all the essential curvatures of order bigger than $1$ of an admissible tuple are trivial. Thus an analytic almost CR manifold has a local CR immersion in $\bC^{\wt m}$ if and only if the total curvature
\beq
\label{stefania}
\Omega^{-1}:\cU\subset M\to Z^{0,2}(\wt\gh_\infty,W)
\eeq
takes values in the space of coboundaries $B^{0,2}(\wt\gh_\infty,W)$. We remark that, by its very definition, the map \eqref{stefania} always takes values in $W\otimes\Lambda^2 W^*$. 
It follows that $(M,\cD,J)$ is locally immersible in $\bC^{\wt m}$ if and only if \eqref{stefania} takes values in $B^{0,2}(\wt\gh_\infty,W)\cap W\otimes\Lambda^2 W^*$ and, by equation \eqref{thee},
if and only if 
\begin{empheq}{align*}
&[X,Y]-[JX,JY]\in\cD\;\;\;\;\;\text{and}\\
&[X,Y]-[JX,JY]=-J[JX,Y]-J[X,JY]\,,
\end{empheq}
for all sections $X,Y$ of $\cD$, i.e. when the almost CR structure is {\it integrable} in the classical sense. We therefore reobtained the classical result by Andreotti and Hill (\cite{AH}) on immersibility of almost CR structures.
\subsection{\it CR-quaternionic  immersions into $\bH$\rm{P}$^n$}\label{CRQ}
Let $(\wt M,Q)$ be an almost quaternionic manifold, a real manifold of dimension $4\wt n$ with  a subbundle $Q$ of $\End(T\wt M)$ locally generated by a triple $(J_1,J_2,J_3)$ of endomorphisms satisfying
$$
J_3=J_1 J_2\;,\;\; J_\alpha^2=-I\;,\;\; J_\alpha\circ J_\beta=-J_\beta\circ J_\alpha\;\;\;\;\text{for all}\;\; 1\leq \alpha,\beta\leq 3\;,\;\; \alpha\neq\beta\,.
$$
To any $(\wt M, Q)$ one can naturally associate a $\wt G$-structure $\wt\pi:\wt P\to \wt M$ with $\wt G=\GL_{\wt n}(\bH)\cdot \Sp_1$ (see \cite{Sal, AM}).
\medskip\par
Let now $M$ be a real manifold of dimension $4\wt n-1$ endowed with an embedding $TM\hookrightarrow E$ of its tangent bundle into a vector bundle $E\to M$ 
of $\rank(E)=4\wt n$ with a linear quaternionic structure on each fibre.
Following \cite{MOP} the pair $(M,E)$ is an {\it almost CR quaternionic manifold (of hypersurface type)}. To any such manifold one can naturally associate a
$G$-structure $\pi:P\to M$ with a particular structure group $G$ which is not semisimple and has Levi factor isomorphic with $\SL_{\wt n-1}(\bH)\times \Sp_1$. We refer to \cite{ASQ} for its explicit description.
\medskip\par
Any hypersurface $M$ of an almost quaternionic manifold $\wt M$ is $\wt P$-regular and inherits an almost CR quaternionic structure, the natural inclusion $TM\subset E:=T\wt M|_M$. Moreover the above mentioned bundle $\pi:P\to M$ coincides with the $G$-structure induced on $M\subset \wt M$ by $\wt\pi:\wt P\to\wt M$. 
\medskip\par
It follows that, given an almost CR quaternionic manifold $(M,E)$ of dimension $4\wt n-1$ and an almost quaternionic manifold $(\wt M,Q)$ of dimension $4\wt n$, the existence of a local CR quaternionic immersion of $(M,E)$ into $(\wt M,Q)$ is equivalent to the local immersibility of $\pi:P\rightarrow M$ into $\wt\pi:\wt P\rightarrow\wt M$.
In particular the existence of a local CR quaternionic immersion of $(M,E)$ in the quaternionic projective space $\wt M=\bH$P$^{\wt n}$ with its standard quaternionic structure $Q=Q_o$ is equivalent to the local immersibility of $\pi:P\rightarrow M$ in the $\wt G$-structure $\wt\pi:\wt P_{o}\to \bH$P$^{\wt n}$ associated with $(\bH$P$^{\wt n},Q_o)$. 
This structure is flat $k$-reductive with $k=2$. We refer to \cite{ASQ} for the explicit expression of the maximal prolongation $\displaystyle\wt\gh_\infty$ of $\wt\gg=\ggl_{\wt n}(\bH)+\sp_1$ and its main properties. In that paper, it is also proved that the groups $H^{1,1}(\wt\gh_\infty,W)$ and $H^{2,1}(\wt\gh_\infty,W)$ are trival. This fact, a careful analysis of the groups $H^{p,2}(\wt\gh_\infty,W)$ and Theorems \ref{topolino}, \ref{r-1} and \ref{quaternion} yield the following result on CR quaternionic immersions. In the statement 
\begin{itemize}
\item[-] $\tE=\bC^{2\wt n-2}$ and $\tH=\bC^{2}$ are the standard representations of $\sl_{2\wt n-2}(\bC)$ and $\sl_{2}(\bC)$ respectively,
\item[-] $\Ad$ is the adjoint representation of $\sl_{2\wt n-2}(\bC)$,
\item[-] $\tD$ is the $\sl_{2\wt n-2}(\bC)$-irreducible module given by the kernel of the natural contraction $\tE\otimes \Lambda^2\tE^*\rightarrow \tE^*$.
\end{itemize}
\begin{theorem}\label{music}
Let $M$ be an almost CR quaternionic manifold of dimension $4\wt n-1\geq 7$ and $\pi:P\rightarrow M$ its canonically associated $G$-structure. Then:
\vskip0.3cm\par\noindent
$i)$ there are natural isomorphisms of $\gs$-modules, $\gs\simeq \sl_{2\wt n-2}(\bC)+\sl_{2}(\bC)$, 
\begin{align*}
H^{0,2}(\wt\gh_\infty,W)\otimes\bC&\simeq \Lambda^2\tE^* S^2\tH+(\tD+\tE^*)S^3\tH+(\Ad+\Lambda^2\tE^*) S^4\tH+\tE^* S^5\tH\,\\
H^{1,2}(\wt\gh_\infty,W)\otimes\bC&\simeq H^{1,2}(\wt\gg)\otimes\bC\,,
\end{align*}
while the cohomology group $H^{2,2}(\wt\gh_\infty,W)$ vanishes;
\vskip0.2cm\par\noindent
$ii)$ there is a canonical $G$-equivariant map $\cR^{1}:P\rightarrow H^{0,2}(\wt\gh_\infty, W)$ which vanishes if $M$ is locally immersible into a quaternionic manifold (i.e. an almost quaternionic manifold with a compatible torsion-free connection);
\vskip0.2cm\par\noindent
$iii)$ if $\cR^1=0$ there is a canonical $G$-equivariant map $\cR^{2}:P\rightarrow H^{1,2}(\wt\gh_\infty,W)$;
\vskip0.2cm\par\noindent
$iv)$  $M$ is locally immersible in $\bH${\rm P}$^n$ around any point if and only if $\cR^1=\cR^2=0$.
\end{theorem}  
\begin{remark}
The fact that the essential curvatures $\cR^1$ and $\cR^2$ of the theorem are intrinsically defined objects relies on the non-trivial fact that $H^{1,1}(\wt\gh_\infty,W)$ and $H^{2,1}(\wt\gh_\infty,W)$ vanish. In general these groups do not vanish: consider the case of isometric immersions into Euclidean spaces where
$$H^{1,1}(\so_{\wt n}(\bR),W)\simeq \underbrace{H^{1,1}(\so_n(\bR),W)}_{=\so_n(\bR)^1=0}+W^\perp\otimes S^2 W^*+\so(W^\perp)\otimes W^*\ .$$
The non vanishing of the last two subspaces correspond to the fact that the existence of a local Riemannian immersion depends on the existence of an appropriate fundamental form and normal connection.
\end{remark}

\end{document}